\documentclass[12pt, twoside, reqno]{amsart}
\calclayout
\usepackage[utf8]{inputenc}
\usepackage{amsmath,amssymb,amscd,amsxtra,latexsym,bbm,graphicx,epsfig,epic,
eepic,oldgerm,bm,psfrag,amsthm,mathrsfs,bm,bbm}
\usepackage{caption}
\usepackage{subcaption}
\usepackage{tikz-cd}
\usepackage{enumitem}
\usepackage[pdftex]{hyperref}

\input xy
\xyoption{all}
\CompileMatrices
    
\theoremstyle{plain}
\newtheorem{theorem}{Theorem}[section]
\newtheorem*{theorem*}{Theorem}
\newtheorem{lemma}[theorem]{Lemma}
\newtheorem{proposition}[theorem]{Proposition}
\newtheorem{corollary}[theorem]{Corollary}
\newtheorem{definition}[theorem]{Definition}
\newtheorem*{definition*}{Definition}

\newtheorem{question}[theorem]{Question}

\newtheorem{assumption}[theorem]{Assumption}
\newtheorem{conjecture}[theorem]{Conjecture}

\theoremstyle{definition}
\newtheorem*{remark*}{Remark}
\newtheorem*{remarks*}{Remarks}
\newtheorem{remark}[theorem]{Remark}

\newtheorem{example}[theorem]{Example}

\newtheorem*{example*}{Example}
\newtheorem*{examples*}{Examples}

\theoremstyle{plain}
\newtheorem{TheoremA}{Theorem}

\newcommand{\proofend}{\hspace*{\fill} $\Box$\\}

\newcommand{\ign}[1]{}

\def\1{\:\!}
\def\2{\;\!}

\def\im{\operatorname {im}}

\def\Diffc0{\operatorname{Diff^c_0}}
\def\Symp{\operatorname{Symp}}

\def\Sympc0{\operatorname{Symp^c_0}}

\def\Int{\operatorname{int}}
\def\Ham{\operatorname{Ham}}
\def\rank{\operatorname{rank}}

\def\GL{\operatorname{GL}}

\def\Lie{\operatorname{Lie}}

\def\Aut{\operatorname{Aut}}
\def\Span{\operatorname{span}}

\def\intl{\operatorname{\ell_{int}}}

\def\CP{\operatorname{CP}}

\def\cd{{\mathcal D}}
\def\ce{{\mathcal E}}

\def\ch{{\mathcal H}}

\def\cj{{\mathcal J}}

\def\cR{{\mathcal R}}

\def\fH{{\mathfrak{H}}}

\def\ft{{\mathfrak{t}}}

\def\fm{{\mathfrak{m}}}

\def\fD{{\mathfrak{D}}}

\def\xh{\widehat{x}}
\def\yh{\widehat{y}}

\def\C{\mathbb{C}}

\def\N{\mathbb{N}}

\def\R{\mathbb{R}}

\def\Z{\mathbb{Z}}

\def\CP{\C P}

\def\U{\operatorname{U}}

\def\pp{\partial}

\def\ddt0{\left. \frac{d}{dt} \right\vert_{t=0}}
\def\dds0{\left. \frac{d}{ds} \right\vert_{s=0}}
\def\ddt{\frac{d}{dt} }
\def\dds{\frac{d}{ds} }
\def\inv{^{-1}}

\def\ni{\noindent}

\def\.{\mskip1mu}
\def\?{\mskip-1mu}

\def\id{\operatorname{id}}

\def\proof{\noindent {\it Proof. \;}}
\newcommand{\proofof}[1]{\ni {\it Proof of #1. }}

\begin{document}

\title[]{Hamiltonian classification of toric fibres and symmetric probes}

\author{Jo\'e Brendel}  
\thanks{Partially supported by the Israel Science Foundation grant 1102/20 and by the ERC Starting Grant 757585}
\address{Jo\'e Brendel,
School of Mathematical Sciences, 
Tel Aviv University
Ramat Aviv, 
Tel Aviv 69978
Israel }
\email{joe@brendel.lu}

\keywords{symplectic geometry, Lagrangian submanifold, toric geometry, toric fibres, probes}

\date{\today}
\thanks{2020 {\it Mathematics Subject Classification.}
Primary 53D12, Secondary 53D20}

\begin{abstract} 
In a toric symplectic manifold, regular fibres of the moment map are Lagrangian tori which are called \emph{toric fibres}. We discuss the question which two toric fibres are equivalent up to a Hamiltonian diffeomorphism of the ambient space. On the construction side of this question, we introduce a new method of constructing equivalences of toric fibres by using a symmetric version of McDuff's probes~\cite{AbrBorMcD14, McD11}. On the other hand, we derive some obstructions to such equivalence by using Chekanov's classification of product tori~\cite{Che96} together with a lifting trick from toric geometry. Furthermore, we conjecture that (iterated) symmetric probes yield all possible equivalences and prove this conjecture for~$\C^n,\CP^2, \C \times S^2, \C^2 \times T^*S^1, T^*S^1 \times S^2$ and monotone~$S^2 \times S^2$. The case of non-monotone~$S^2 \times S^2$ will be treated in~\cite{BreKim23}.

This problem is intimately related to determining the \emph{Hamiltonian monodromy group} of toric fibres, i.e.\ determining which automorphisms of the homology of the toric fibre can be realized by a Hamiltonian diffeomorphism mapping the toric fibre in question to itself. For the above list of examples, we determine the Hamiltonian monodromy group for all toric fibres.
\end{abstract}

\maketitle

\section{Introduction}

\subsection{Symmetric probes}
Probes were introduced by McDuff~\cite{McD11} to prove that some toric fibres are displaceable. Probes are rational segments in the base of a toric base polytope which hit the boundary integrally transversely in one point. The latter condition implies that one can perform symplectic reduction on a probe and obtain a two-disk as reduced space, see also the exposition in~\cite{AbrMac13}. Toric fibres map to circles in the reduced space, where displaceability questions are easy to settle since they boil down to area arguments. 

\begin{figure}
	\begin{center}
	\begin{tikzpicture}
		\node[inner sep=0pt] at (0,0)
    			{\includegraphics[trim={1.5cm 6.5cm 4cm 2cm},clip,scale=0.75]						{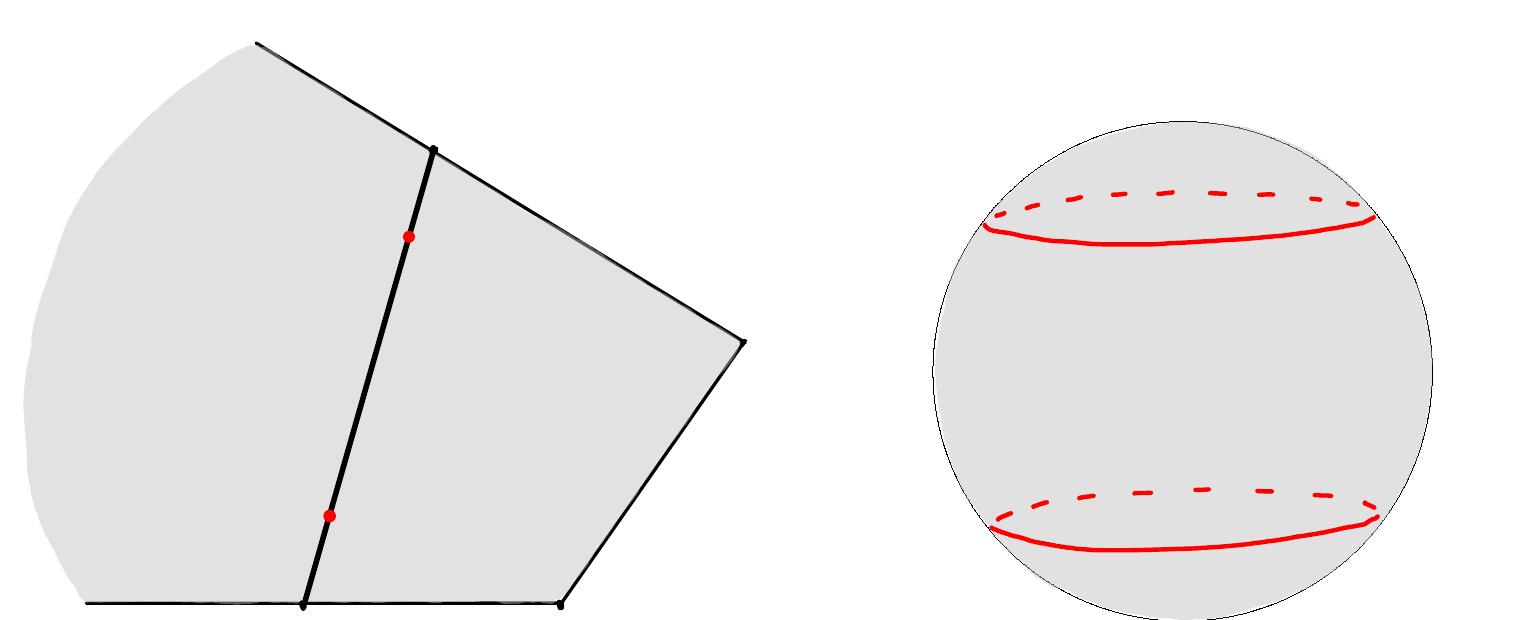}};
    		\node at (-2.4,0.55){$x$};
    		\node at (-2.9,-1.5){$x'$};
    		\node at (-4,1.3){$\Delta$};
    		\node at (4.5,1.5){$S^2$};
    		\node at (-2.7,-0.5){$\sigma$};
    		\node at (1.3,0.64){$S^1_x$};
    		\node at (1.35,-1.675){$S^1_{x'}$};
	\end{tikzpicture}
	\caption{A symmetric probe~$\sigma$ in a moment polytope~$\Delta$ with points~$x,x'$ at equal distance to the boundary. The toric fibres~$T(x),T(x')$ map to the circles~$S^1_x,S^1_{x'} \subset S^2$ under symplectic reduction.}
	\label{fig:1}
	\end{center}
\end{figure}
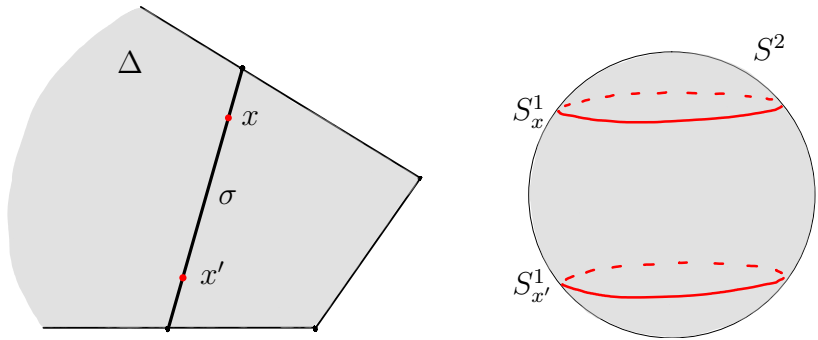

Symmetric probes are rational segments in which \emph{both} endpoints hit the boundary of the moment polytope integrally transversely. They were introduced in a follow-up paper to~\cite{McD11} by Abreu--Borman--McDuff~\cite{AbrBorMcD14} to settle some more subtle displaceability questions. Here, we use them to a different end. The reduced space associated to a symmetric probe is a two-sphere and the quotient map takes toric fibres to orbits of the standard~$S^1$-action on the two-sphere. Observe that -- except for the equator -- orbits of this circle action in~$S^2$ appear in pairs which are Hamiltonian isotopic. Our main observation is that, since Hamiltonian isotopies in reduced spaces can be lifted, this proves that toric fibres corresponding to such pairs of circles are Hamiltonian isotopic, as well. This is illustrated in Figure~\ref{fig:1}.

To state this formally, let us introduce some notation. Let~$(X^{2n},\omega)$ be a (not necessarily compact) toric symplectic manifold with moment map~$\mu \colon X \rightarrow \R^n$ and moment polytope~$\mu(X) = \Delta$. For~$x \in \Int \Delta$ the set~$T(x) = \mu^{-1}(x)$ is a Lagrangian torus, called a \emph{toric fibre}. A \emph{symmetric probe}~$\sigma \subset \Delta$ is a rational segment intersecting~$\partial \Delta$ integrally transversely in the in interior of two facets, see also~\cite[Definition 2.2.3]{AbrBorMcD14}. An intersection of a rational line and a rational hyperplane is called \emph{integrally transverse} if their union contains a~$\Z$-basis of~$\Z^n$. See also Definition~\ref{def:symmprobe} and the surrounding discussion or~\cite[\S 2.1]{McD11} for more details. 

\begin{TheoremA}
\label{thm:mainA}
Let~$(X,\omega)$ be a toric symplectic manifold and let~$\sigma \subset \Delta$ be a symmetric probe in its moment polytope. Furthermore, let~$x,x' \in \sigma$ be two points at equal distance to the boundary~$\pp \Delta$. Then~$T(x)$ and~$T(x')$ are Hamiltonian isotopic. 
\end{TheoremA}

\subsection{Classification of toric fibres}
\label{ssec:introclass}
Deciding which two given Lagrangians~$L$ and~$L'$ in $(X,\omega)$ can be mapped to one another by a symplectomorphism or by a Hamiltonian diffeormorphism is a central question in symplectic geometry. In many situations, it is quite hopeless to give a full classification -- even constructing examples of Lagrangians that are not equivalent to known ones (so-called \emph{exotic} Lagrangians) is an active area of research where many questions are open, see for example~\cite{Aur15, CheSch10, Via17, Via16}. In this paper we care about the following classification question of Lagrangian submanifolds.

\begin{question}
\label{qu:class}
In a toric symplectic manifold~$(X,\omega)$, give a classification of toric fibres up to Hamiltonian diffeomorphisms of the ambient space.
\end{question}

\begin{remark}
One can ask the same questions for symplectomorphisms of the ambient space. In this paper we focus on the case of Hamiltonian diffeomorphisms. See also Remark~\ref{rk:sympclass}.
\end{remark}

Although Question~\ref{qu:class} is a much less ambitious question than a full classification of all Lagrangian tori (since we exclude exotic tori a priori) of~$X$, it is open except for a few special cases and surprisingly absent from the literature. To our knowledge, it has only been answered for~$\C^n$ (where toric fibres are simply product tori) by Chekanov~\cite[Theorem A]{Che96} and for~$\CP^2$ by Shelukhin--Tonkonog--Vianna~\cite[Proposition 7.1]{SheTonVia19}. 

Let us make some conventions. From now on, we call~$T(x)$ and~$T(x')$ \emph{equivalent} and write~$T(x) \cong T(x')$ if they can be mapped to one another by a Hamiltonian diffeomorphism of the ambient space. Furthermore, let 
\begin{equation}
	\fH_x = \{ x' \in \Int \Delta \, \vert \, T(x) \cong T(x') \},
\end{equation}
the set of toric fibres equivalent to~$T(x)$. A first guess may be that~$\fH_x = \{x\}$, since the zero-section in~$T^*T^n$ is non-displaceable (see~\cite[\S 11.3]{McDSal17} and the references therein) and thus this is true if we restrict our attention to Hamiltonian diffeomorphisms supported in a Weinstein neighbourhood of~$T(x)$. However, a glance at~$S^2$ shows that this guess is wrong, since one can use the topology of the ambient space to obtain non-trivial equivalences of toric fibres. More generally, by Theorem~\ref{thm:mainA}, symmetric probes (and their concatenations) can be used to construct equivalences of toric fibres up to Hamiltonian diffeomorphisms. Let us also point out that symmetric probes are abundant in arbitrary toric manifolds -- at least close to the boundary of the moment polytope, see \S\ref{ssec:arbitrary}. We conjecture that the method of constructing equivalent toric fibres by symmetric probes gives a complete anwer to the classification question.

\begin{conjecture}
\label{conj:main}
Two toric fibres~$T(x),T(x') \subset X$ are equivalent if and only if they are equivalent by a sequence of symmetric probes. 
\end{conjecture}

In Section~\ref{sec:examples}, we verify this conjecture for~$\C^n$ and~$\CP^2$ (where the classification was previously known), for~$\C \times S^2, \C^2 \times T^*S^1, T^*S^1 \times S^2$ and for monotone~$S^2 \times S^2$ (where we classify toric fibres). The classification of toric fibres in non-monotone~$S^2 \times S^2$ is more intricate and is given in~\cite{BreKim23}.  

On the side of obstructions to Hamiltonian equivalence, we prove the following. 

\begin{TheoremA}
\label{thm:mainB}
If toric fibres~$T(x),T(x') \subset X$ of a compact toric manifold~$X$ are Hamiltonian isotopic, then the following three invariants agree
\begin{equation}
\label{eq:chekinvariants}
d(x) = d(x'), \quad 
\#_d(x) = \#_d(x'), \quad 
\Gamma(x) = \Gamma(x').
\end{equation}
\end{TheoremA}

The invariant~$d(x) \in \R$ is the \emph{integral affine distance} of~$x$ to the boundary of the moment polytope. The invariant~$\#_d(x) \in \N_{\geqslant 1}$ is the number of facets of~$\Delta$ realizing the minimal distance~$d(x)$. Both of these invariants are \emph{hard} in the symplectic sense. The last invariant is the subgroup
\begin{equation}
	\Gamma(x) = \Z\langle \ell_1(x)-d(x), \ldots, \ell_N(x)-d(x) \rangle \subset \R
\end{equation}
and it is soft. Here~$\ell_i(x)$ denotes the integral affine distance of~$x$ to the~$i$-th facet of~$\Delta$. Since these invariants are derived from Chekanov's invariants~\cite[Theorem A]{Che96} of product tori in~$\R^{2n} = \C^n$, we call them \emph{Chekanov invariants}. 

Let us outline the proof of Theorem~\ref{thm:mainB}. Suppose~$T(x), T(x') \subset X$ are Hamiltonian isotopic fibres. By a construction going back to Delzant~\cite{Del88}, we can view~$X$ as a symplectic quotient of~$\C^N$, where~$N$ is the number of facets of~$\Delta$. The preimages of the tori~$T(x),T(x')$ under the symplectic quotient map are the product tori~$T(\ell(x)), T(\ell(x')) \subset \C^N$, where~$\ell = (\ell_1, \ldots, \ell_N)$. The Hamiltonian isotopy mapping~$T(x)$ to~$T(x')$ lifts to a Hamiltonian isotopy of~$\C^N$ mapping~$T(\ell(x))$ to~$T(\ell(x'))$. This means that Chekanov's invariants for product tori have to agree on~$T(\ell(x))$ and~$T(\ell(x'))$, which yields the statement. To our knowledge, this \emph{lifting trick} first appeared in~\cite{AbrMac13} to prove non-dispaceability of certain fibres and it was also heavily used in~\cite{Bre20}. It is not obvious to us how to prove Theorem~\ref{thm:mainB} directly, i.e.\ without using the lifting trick. The first two invariants are clearly related to the area and the number of non-trivial Maslov two $J$-holomorphic disks of minimal area with boundary on the corresponding tori, respectively. It is not obvious how to pursue this due to the lack of monotonicity, although an approach in the spirit of~\cite{SheTonVia19} may be promising, especially in dimension four, see Remark~\cite[\S 5.6]{SheTonVia19}. 

The invariants in Theorem~\ref{thm:mainB} are not complete, even in very simple examples such as~$\CP^2$, see Example~\ref{ex:cp2notcomplete}. We suspect that the first two invariants are all there is in terms of hard obstructions, but that the soft invariant~$\Gamma(\cdot)$ is far from optimal -- this is the case in all examples where we know the classification.

\subsection{Examples}
\label{ssec:introex}

Let us give some examples of symmetric probes. In dimension two, there are not many toric spaces. The main examples are~$T^*S^1 = S^1 \times \R$ equipped with the standard exact symplectic form and moment map given by projection to the~$\R$-coordinate;~$\C = \R^2$ equipped with the standard symplectic form and moment map~$z \mapsto \pi \vert z \vert^2$, and~$S^2$ equipped with the height function. We normalize the symplectic form~$\omega_{S^2}$ such that~$\int_{S^2} \omega_{S^2} = 2$ meaning that the corresponding moment polytope is~$[-1,1]$. In the two-dimensional setting, symmetric probes are not interesting and the classification of toric fibres boils down to simple area arguments. However, some four-dimensional products of the above examples (equipped with the product symplectic and toric structures) already contain non-trivial probes. 

\begin{figure}
	\begin{center}
	\begin{tikzpicture}
		\node[inner sep=0pt] at (0,0)
    			{\includegraphics[trim={3cm 6.5cm 1cm 3cm},clip,scale=0.75]						{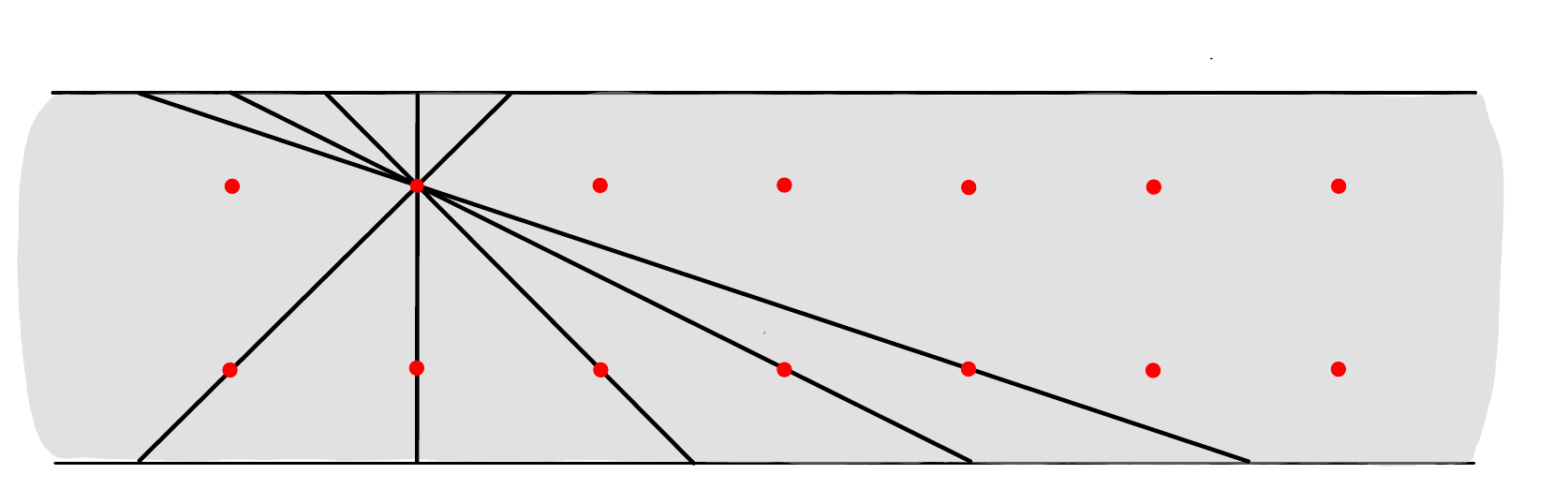}};
    		\node at (-2.2,0.6){$(x,y)$};
	\end{tikzpicture}
	\caption{The set~$\fH_{(x,y)}$ for~$T(x,y) \subset T^*S^1 \times S^2$ and some symmetric probes.}
	\label{fig:2}
	\end{center}
\end{figure}

In~$\C^2$, there is one non-trivial probe in direction~$(1,-1)$, which can be used to show that~$T(x,y) \cong T(y,x)$, which also follows from the fact that all elements in~$\U(2)$ can be realized by Hamiltonian diffeomorphisms. These are all possible equivalences in~$\C^2$, as was shown by Chekanov~\cite{Che96}. In~$T^*S^1 \times S^2$, all directions~$(k,1)$ for~$k \in \Z$ give symmetric probes, see Figure~\ref{fig:2}. This proves that~$T(x,y)$ is Hamiltonian isotopic to all~$T(x+2ky,\pm y)$. Note that this also follows from a suspension argument due to Polterovich, see~\cite[Example 6.3.C]{Pol01} and the discussion in~\cite[\S1.3]{MakSmi19}. The example~$\C \times S^2$ is obtained from the previous one by a vertical symplectic cut and we will see in~\S\ref{ssec:excs2} that there are slightly more equivalences between toric fibres. In~$\CP^2$ and monotone~$S^2 \times S^2$, it is easy to see that symmetric probes realize all equivalences of toric fibres coming from symmetries of the moment polytope. In all of these examples, the method by probes is sharp and the classification of toric fibres is discussed in detail in Section~\ref{sec:examples}. 

In dimensions~$\geqslant 6$, the situation is quantitatively different from the above examples. Indeed, the set~$\fH_x$ has accumulation points in~$\Delta$ for many~$x \in \Int \Delta$, see Corollary~\ref{cor:accumulation}. This already occurs in the case of~$\C^3$, treated by Chekanov~\cite{Che96}, see also Theorem~\ref{thm:chekanov}. In essence, this is due to the existence of a symmetric probe in direction~$(1,1,-1)$ (or coordinate permutations thereof), see Figure~\ref{fig:5}. In~\S\ref{ssec:chekrev}, we show that one can recover Chekanov's classification using symmetric probes. The property that~$\fH_x$ has accumulation points is not exclusive to dimension six and above. In fact, in the forthcoming~\cite{BreKim23}, we show that this occurs in~$S^2 \times S^2$ equipped with any non-monotone symplectic form. 

In light of this, it would be very interesting to characterize the toric manifolds having the property that there exists~$x \in \Int \Delta$ with~$\fH_x$ not discrete.

\subsection{Hamiltonian monodromy of toric fibres} Let~$x \in \sigma \subset \Delta$ be the midpoint of a symmetric probe. The corresponding toric fibre~$T(x)$ projects to the equator of the sphere obtained as a reduced space and thus we do not get any equivalence with another toric fibre by the above method. However, we still get information about~$T(x)$. Indeed, we can lift a Hamiltonian isotopy mapping the equator in the reduced sphere to itself but changing the orientation of the equator. By lifting the such a Hamiltonian isotopy, we obtain a Hamiltonian isotopy mapping~$T(x)$ to itself with non-trivial homological monodromy, meaning that it induces a non-trivial map in~$\Aut H_1(T(x);\Z)$. An explicit formula for this monodromy map in terms of data related to the symmetric probe~$\sigma$ is given in~\eqref{eq:basicinvolution}. 

\begin{definition}
\label{def:hmg}
Let~$L \subset (X,\omega)$ be a compact Lagrangian submanifold. The \emph{Hamiltonian monodromy group} is given by
\begin{equation}
\label{eq:hmg}
	\ch_L = \{ (\phi\vert_L)_* \in \Aut H_1(L;\Z) \; \vert \; \phi \in \Ham(X,\omega),\, \phi(L) = L \}.
\end{equation}
\end{definition}

The analogous monodromy group for symplectomorphisms was computed by Chekanov for product tori and Chekanov tori in~\cite[Theorem 4.5]{Che96} and, in that case, the Hamiltonian monodromy group actually agrees with it. To our knowledge this is the first occurence of this kind of question in the literature. See also Yau~\cite{Yau09} for related results and Hu--Lalonde--Leclercq~\cite{HuLalLec11} which establishes that weakly exact Lagrangian manifolds have trivial Hamiltonian monodromy group. See Porcelli~\cite{Por22} for recent progress in the same direction. Another recent work is Augustynowicz--Smith--Wornbard~\cite{AugSmiWor22} which makes significant progress in case~$L$ is a monotone Lagrangian torus and provides an excellent overview of the topic in its introduction. 

Let~$\xi_1,\ldots,\xi_N \in \Z^n$ be the set of inward pointing primitive normal vectors to the facets of~$\Delta$, as in~\eqref{eq:polytope}. The vectors~$\xi_i$ naturally determine homology classes~$\xi_i \in H_1(T(x))$ for every toric fibre~$T(x)$. See for example the discussion surrounding~\eqref{eq:lambda}. Let~$\cd(x)$ be the subset of those normal vectors realizing the minimal integral affine distance of~$x$ to facets,
\begin{equation}
	\cd(x) = \{\xi_i \, \vert \, \ell_i(x) = d(x) \}.
\end{equation}
We call elements of this subset \emph{distinguished classes}. Note that~$\# \cd(x) = \#_d(x)$. The following is an obstruction result for Hamiltonian monodromy of toric fibres.

\begin{TheoremA}
\label{thm:mainC}
Let~$T(x) \subset X$ be a toric fibre in a compact toric manifold. Every element in the Hamiltonian monodromy group~$\ch_{T(x)}$ acts by a permutation on the set~$\cd(x)$ of distinguished classes.
\end{TheoremA}

This theorem again follows from Chekanov's work~\cite[Theorem 4.5]{Che96} together with the lifting trick discussed in~\S\ref{ssec:introclass}. In fact, we get a stronger statement, see Theorem~\ref{thm:ambientmonodromy}. The number~$\#_d(x)$ of distinguished classes is maximal if~$T(x)$ is monotone, since all integral affine distances are equal in that case. In fact, in the monotone case, we recover~\cite[Theorem 2]{AugSmiWor22} for Hamiltonian diffeomorphisms, see Corollary~\ref{cor:smith}. Note that Theorem~\ref{thm:mainC} does not require monotonicity. 

In terms of examples, we give a complete description of~$\ch_L$ for all toric fibres in~$S^2 \times S^2, \CP^2, \C \times S^2, \C^2 \times T^*S^1$ and~$T^*S^1 \times S^2$, and show that all Hamiltonian monodromy elements can be realized by symmetric probes as outlined above.

\subsection{Outline}
In Section~\ref{sec:toric}, we review the relevant toric geometry and in particular we discuss \emph{toric reduction}, a version of symplectic reduction which is compatible with the toric structure and on which we rely to prove Theorems~\ref{thm:mainA},~\ref{thm:mainB} and~\ref{thm:mainC}. Section~\ref{sec:symmprobes} is the heart of this paper, where we discuss symmetric probes and prove Theorem~\ref{thm:mainA}. In Section~\ref{sec:chekanovinv}, we discuss obstructions to the equivalence of toric fibres and prove Theorem~\ref{thm:mainB}. Furthermore, we discuss obstructions to which Hamiltonian monodromy can be realized for toric fibres. Section~\ref{sec:examples} is dedicated to examples and serves to illustrate the results of the previous sections. 

\subsection*{Acknowledgements}
We cordially thank Jonny Evans, Joontae Kim and Felix Schlenk for many useful discussions. We are grateful to Jack Smith for generously sharing his insights on Hamiltonian monodromy and offering explanations about~\cite{AugSmiWor22}. This work was started at Université de Neuchâtel, partially supported by SNF grant 200020-144432/1, and continued at Tel Aviv University, partially supported by the Israel Science Foundation grant 1102/20 and by the ERC Starting Grant 757585.

\section{Some toric symplectic geometry}
\label{sec:toric}

In this section, we review toric geometry with special emphasis on a certain type of symplectic reduction, which we call \emph{toric reduction}. Toric reduction generalizes probes as well as Delzant's construction of toric symplectic manifolds, both of which heavily feature in this paper.

\subsection{Toric manifolds}
A symplectic manifold~$(X^{2n},\omega)$ together with a moment map~$\mu \colon X \rightarrow \ft^*$ is called \emph{toric} if~$\mu$ generates an effective Hamiltonian action of the~$n$-torus~$T^n$. By~$\ft^*$ we denote the dual of the Lie algebra~$\ft$ of~$T^n$. Choosing an identification~$T^n \cong \R^n/\Z^n$ induces an identification~$\ft^* \cong \R^n$ and, depending on context, we will use both the invariant way and the coordinate-dependent way of seeing things. Note that some symplectic manifolds admit distinct toric structures and hence we are really concerned with the triple~$(X,\omega,\mu)$ when we say \emph{toric manifold} although we may just write~$X$ or~$(X,\omega)$ for simplicity. 

A classical result by Delzant~\cite{Del88} states that if~$X$ is \emph{compact} toric\footnote{Many authors include compactness in the definition of \emph{toric}, but we do not.}, then the image~$\Delta = \mu(X)$ is a so-called \emph{Delzant polytope}, and that Delzant polytopes (up to integral affine transformations) classify toric manifolds up to equivariant symplectomorphism. There are many classical references for toric manifolds, e.g.~\cite{Aud04, Can03, Gui94}, and we refer to these for details. We revisit part of Delzant's result in~\S\ref{ssec:delzant}. 

Due to Delzant's theorem, the moment polytope associated to a toric manifold~$X$ is a crucial object of study. We view it as
\begin{equation}
	\label{eq:polytope}
	\Delta = \{ x \in \ft^* \;\vert\; \ell_i(x) \geqslant 0 \}, \quad
	\ell_i(x) = \langle x , \xi_i \rangle + \lambda_i.
\end{equation}
Here, we view the vectors~$\xi_i$ in~$\ft$ and~$\langle \cdot, \cdot \rangle$ denotes the natural pairing of~$\ft$ and its dual. Note that~$\ft$ contains a natural lattice~$\Lambda$ obtained as the kernel of the exponential map~$\exp \colon \ft \rightarrow T^n$. Similarly, the dual~$\ft^*$ contains the dual lattice~$\Lambda^*$. If we choose a basis, we can identify~$\Lambda \cong \Z^n$ and dually~$\Lambda^* \cong \Z^n$. Again, depending on context, we use both the invariant viewpoint and the coordinate-dependent one. Furthermore, since~$\Delta$ is rational (with respect to~$\Lambda^*$), we can choose the vectors~$\xi_i$ to be primitive in~$\Lambda$.

\begin{definition}
A vector~$v \in \Lambda$ in a lattice~$\Lambda$ is called \emph{primitive} if~$\alpha v \notin \Lambda$ for all~$0<\alpha<1$.
\end{definition}

Together with~\eqref{eq:polytope}, this condition uniquely determines~$\xi_i$ and~$\lambda_i$ in terms of~$\Delta$ and vice-versa. As we have mentioned above,~$\Delta$ is a \emph{Delzant polytope}, meaning that at every vertex the vectors~$\xi_i$ determining the facets meeting at that vertex form a basis of the lattice~$\Lambda$ over the integers. There is a natural symmetry group acting on~$\Delta \subset \ft^*$ without changing the toric manifold determined by~$\Delta$. 

\begin{definition}
The \emph{integral affine transformations} of~$(\ft^*,\Lambda^*)\cong(\R^n,\Z^n)$ are the elements in the group
\begin{equation}
	\Aut \Lambda^* \ltimes \ft^* \cong \GL(n;\Z) \ltimes \R^n.
\end{equation}
\end{definition}

The elements in~$\Aut \Lambda^* \cong \GL(n;\Z)$ correspond to base changes in the torus~$T^n$, whereas the translation part~$\ft^* \cong \R^n$ correspond to adding constant elements to the moment map. Neither of these transformations changes the Hamiltonian~$T^n$-action.

\subsection{Toric reduction} 
\label{ssec:toricred} 

In this paragraph we are interested in symplectic reduction with respect to subtori of a toric~$T^n$-action. We call symplectic reduction of this type \emph{toric reduction}. The symplectic quotient of this operation inherits a toric structure with moment polytope obtained by intersecting~$\Delta$ with an affine rational subspace in~$\ft^*$. Roughly speaking, toric reductions are in bijection with inclusions (which are compatible in the sense of Defintion~\ref{def:redadm}) of the moment polytope of the reduced space into the moment polytope of the initial space. Although we could not find a precise statement of sufficient generality in the literature, this idea is hardly new - see for example~\cite{AbrMac13}. In fact, as we will discuss in~\S\ref{ssec:delzant}, the Delzant construction and McDuff's probes are special cases of Theorem~\ref{thm:toricred}. What may be new is the precise formulation we give in Definition~\ref{def:redadm} of the conditions for this reduction to yield a smooth symplectic quotient in terms of the geometry of~$\Delta$.

Let~$X$ be a toric manifold and~$\Delta$ its moment polytope. Note that symplectic reduction with respect to the full~$T^n$-action is pointless. Indeed, the reduced spaces are zero-dimensional. However, it is quite fruitful to perform symplectic reduction with respect to a subtorus~$K \subset T^n$. Dually, we may look at affine rational subspaces~$V \subset \R^n \cong \ft^*$. Indeed, to any affine rational subspace~$V$ we can associate its complementary torus
\begin{equation}
\label{eq:comptorus}
K_V = \exp(V^0), \quad
V^0 = \{ \xi \in \ft \,\vert\, \langle x - x', \xi \rangle = 0, \; x,x' \in V \} \subset \ft,
\end{equation}
and vice-versa. Rationality of~$V$ is equivalent to the compactness of~$K_V$. The subspace~$V$ is a level set of the natural projection~$\ft^* \rightarrow \Lie(K_V)^*$, meaning that the moment map~$\mu_{K_V} \colon X \rightarrow \Lie(K_V)^*$ generating the induced~$K_V$-action on~$X$ has level set~$\mu\inv(\Delta \cap V)$ for some suitable level. Thinking in terms of~$V$ and instead of~$K_V$ or~$\mu_{K_V}$ has the advantage that both the subtorus and the level at which we wish to carry out reduction are fixed by a choice of~$V$. Furthermore, one can easily read off the integral affine geometry of the pair~$(\Delta,V)$ whether the action of~$K_V$ on~$\mu\inv(\Delta \cap V)$ is free (and hence the reduction admissible). Obviously, this is not always the case, since~$V$ may contain a vertex of~$\Delta$ for example. 

\begin{definition}
\label{def:redadm}
Let~$\Delta$ be a Delzant polytope and let~$V$ be an affine rational subspace. We call the pair~$(\Delta,V)$ \emph{reduction-admissible} if, for every face~$F \subset \Delta$ intersecting~$V$, the union of (the linear part of)~$F$ and (the linear part of)~$V$ contains a basis of the lattice~$\Lambda^*$.
\end{definition}

Analogously we call a polytope~$\Delta' \subset \Delta$ \emph{reduction-admissible} if it is obtained as the intersection~$\Delta' = \Delta \cap V$ of~$\Delta$ with a reduction-admissible~$V$. Note that one only needs to check reduction-admissibility at the faces~$F$ of the smallest dimension for which~$V \cap F$ is non-empty, i.e.\ at the vertices of the polytope~$\Delta'$.

\begin{theorem}[Toric reduction]
\label{thm:toricred}
Let~$\Delta \subset \R^n$ be a Delzant polytope and~$V \subset \R^n$ an affine rational subspace such that the pair~$(\Delta,V)$ is reduction-admissible. Then the action of~$K_V = \exp(V^0)$ on~$Z = \mu\inv(\Delta \cap V)$ is free and the reduced space~$X' = Z/K_V$ is itself toric with moment polytope~$\Delta' = \Delta \cap V$.
\end{theorem}

\proof
Let~$e_1^*,\ldots,e_n^* \in \R^n = \ft^*$ be the standard basis. Reduction-admissibility implies that, up to applying an integral affine transformation, we may assume that 
\begin{equation}
	V = \Span_{\R}\{e_1^*,\ldots,e_i^*\}, \quad
	F = \Span_{\R}\{e_j^*,\ldots,e_n^*\}, \quad j \leqslant i + 1. 
\end{equation}
In this normal form, we have~$V^0 = \Span_{\R}\{e_{i+1},\ldots,e_n\}$ and hence~$K_V = \{1\} \times T^{n-i}$. This subtorus acts freely on~$\mu\inv(F)$. Since this holds for any facet~$F$ intersecting~$V$, the action of~$K_V$ is free and thus symplectic reduction is admissible.

The quotient manifold carries a residual~$T^n/K_V$-action. It is effective, since the~$T^n$-action on~$X$ is. Since~$\mu$ is invariant under the~$T^n$-action, it is in particular invariant under the induced~$K_V$-action and thus its restriction to~$Z = \mu\inv(\Delta \cap V)$ factors through the quotient by~$K_V$ and has image~$\Delta' = \Delta \cap V$. It is not hard to check that the map obtained in this way is a moment map generating the~$T^n/K_V$-action on the quotient. For dimensional reasons, the resulting action is toric.
\proofend

Let~$M = T^n/K_V$ be the torus acting by the residual action. Note that, by definition,~$\Delta'$ is contained in~$\ft^*$ instead of~$\Lie(M)^* = \fm^*$. However, one can pick an identification of~$(\fm^*, \Lambda_{M}^*)$ with~$(V,\Lambda \cap V)$ and, up to an element in the integral affine transformations of~$(\fm^*,\Lambda_M^*)$, this yields a well-defined polytope~$\Delta' \subset \fm^*$. Conversely, given an integral affine embedding
\begin{equation}
	\label{eq:iota}
	\iota \colon (\fm^*,\Lambda_M) \hookrightarrow (\ft^*, \Lambda), \quad \iota(\Delta') = \iota(\ft^*) \cap \Delta
\end{equation}
such that~$(\Delta,\iota(\fm^*))$ is reduction-admissible, there is a symplectic reduction from~$X$ to~$X'$. To summarize, there is a short exact sequence of tori,
\begin{equation}
\label{eq:sestori}
 0 \rightarrow K_V \hookrightarrow T^n \stackrel{\Xi}{\rightarrow} M \rightarrow 0,
\end{equation}
where~$T^n$ acts on~$X$ and~$M$ acts on the reduced space~$X'$ such that the reduction map~$p \colon Z \rightarrow X'$ is equivariant with respect to the~$T^n$- and~$M$-actions meaning that
\begin{equation}
	\label{eq:pequivariance}
	p(t.x)=\Xi(t).p(x), \quad t \in T^N, x \in X.
\end{equation}
In particular, orbits are mapped to orbits under toric reduction. This will be used in~\S\ref{ssec:toricfibres}.

\subsection{Delzant construction}
\label{ssec:delzant}
The Delzant construction gives a recipe for constructing a toric manifold~$(X,\omega,\mu)$ from a compact Delzant polytope~$\Delta$. We review it here, since it will be used in Section~\ref{sec:chekanovinv}, and refer to~\cite{Gui94} for details. Actually, the Delzant construction is a special case of toric reduction as discussed in~\S\ref{ssec:toricred} where~$X$ is obtained as a symplectic quotient of some~$\C^N$ equipped with its standard toric structure. 

Let~$\Delta \subset \ft^*$ be a Delzant polytope with~$N$ facets. Since~$\Delta$ is compact, we have~$N > n$. Let~$(\C^N,\omega_0)$ be the standard symplectic vector space equipped with the moment map
\begin{equation}
	\label{eq:stdmm}
	\mu_0 \colon \C^N \rightarrow (\ft^N)^* \cong \R^N, \quad
	(z_1,\ldots,z_N) \mapsto (\pi\vert z_1 \vert^2, \ldots, \pi\vert z_N \vert^2),
\end{equation}
which generates the standard~$T^N$-action on~$\C^N$ by rotation in the factors. Its image is the positive orthant~$\R^{N}_{\geqslant 0}$. Instead of starting with the subtorus~$K \subset T^N$ by which to reduce, we start by defining an inclusion
\begin{equation}
	\label{eq:ell}
	\ell \colon \ft^* \hookrightarrow \R^N, \quad
	x \mapsto (\ell_1(x),\ldots,\ell_N(x)),
\end{equation}
which maps~$\Delta$ to~$\R^N_{\geqslant 0}$. The components~$\ell_i$ defined in~\eqref{eq:polytope} are the functions measuring the integral affine distance of a given point to the facets of~$\Delta$. The map~$\ell$ is an integral affine embedding as in~\eqref{eq:iota} and the subtorus~$K$ by which we reduce is given~$K = \exp (\im \ell)^0 \subset T^N$. Using the Delzant condition on~$\Delta$, it is easy to check that the inclusion~$\ell(\Delta) \subset \R^N_{\geqslant 0}$ is admissible in the sense of Definition~\ref{def:redadm}. Thus the toric symplectic manifold~$(X,\omega,\mu)$ is obtained as symplectic quotient~$X = \mu_0^{-1}(\ell(\Delta)) / K$. 

Let us illustrate this by a simple example.

\begin{figure}
	\begin{center}
	\begin{tikzpicture}
		\node[inner sep=0pt] at (0,0)
    			{\includegraphics[trim={3cm 5cm 5.5cm 2.5cm},clip,scale=0.75]						{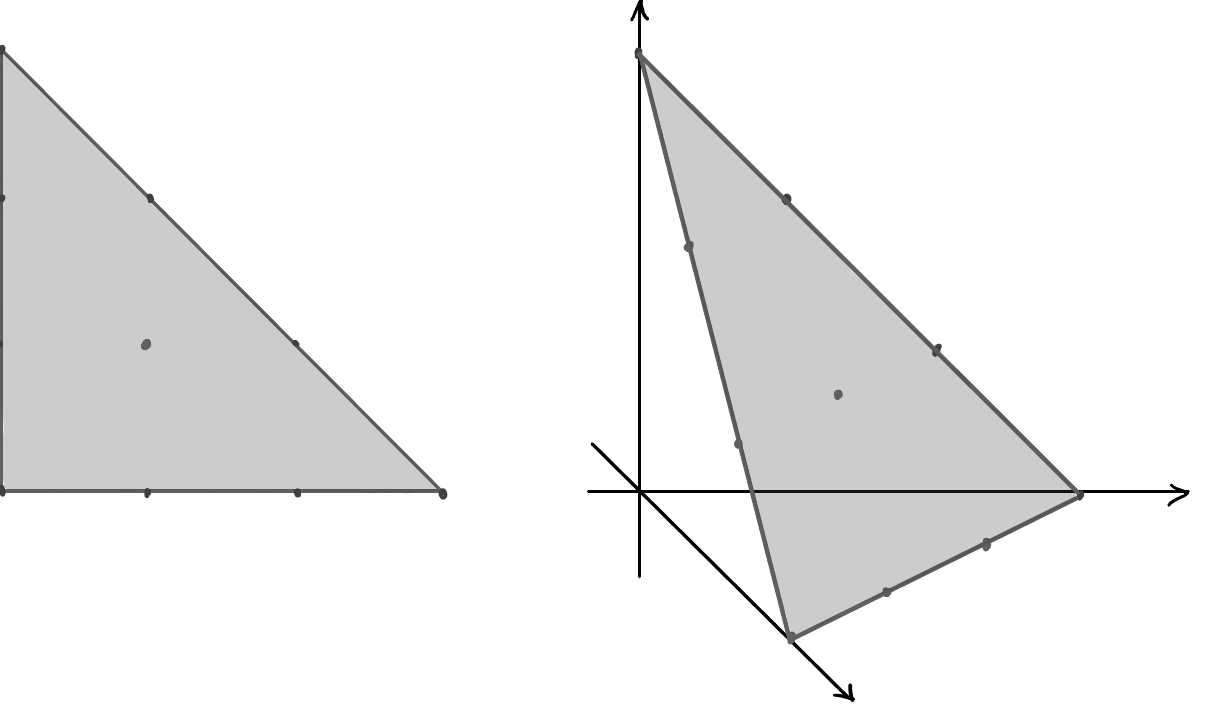}};
    		\node at (-3.8,2){$\Delta $};
    		\node at (1.25,2){$\ell(\Delta)$};
	\end{tikzpicture}
	\caption{The idea of the Delzant construction in the case of~$X = \CP^2$. The complement of~$\im \ell$ generates the circle action by which the symplectic reduction is performed.}
	\label{fig:3}
	\end{center}
\end{figure}

\begin{example}[Complex projective plane]
\label{ex:cp2}
Let~$\Delta \subset \ft^* = \R^2$ be the simplex defined by
\begin{equation}
	\ell_1(x) = x_1 + 1, \quad
	\ell_2(x) = x_2 + 1, \quad
	\ell_3(x) = - x_1 - x_2 +1.
\end{equation}
This simplex is Delzant and since~$N=3$, we will obtain~$X$ as a symplectic reduced space of~$\C^3$. The map~$\ell$ is depicted in Figure~\ref{fig:3}. The orthogonal complement~$(\im \ell)^{\perp}$ is spanned by~$(1,1,1)$ and thus~$K=\{(t,t,t) \,\vert\, t \in S^1\} \subset T^3$ and~$\mu_K(z) = \pi(\vert z_1 \vert^2 + \vert z_2 \vert^2 + \vert z_3 \vert^2)$. We conclude that the symplectic reduction~$\mu_K\inv(3) = S^5(3) \rightarrow \CP^2$ corresponds to the Hopf fibration map. The symplectic form one obtains by this procedure is the Fubini--Study form~$\omega_{\CP^2}$ with normalization~$\int_{\CP^1}\omega_{\CP^2} = 3$.
\end{example} 

\subsection{Toric fibres} 
\label{ssec:toricfibres} 
Every toric manifold~$X^{2n}$ contains an~$n$-parametric family of Lagrangian tori called \emph{toric fibres}.

\begin{definition}
Let~$x \in \Int\Delta$ be a point in the interior of a toric moment polytope. The corresponding preimage~$T(x)=\mu\inv(x)$ is called a \emph{toric fibre}. 
\end{definition}

\begin{example}
\label{ex:prodtori}
{\rm (Product tori)
The toric fibres of the standard toric structure~\eqref{eq:stdmm} are \emph{product tori}~$\mu_0\inv(a_1,\ldots,a_N) = S^1(a_1)\times \ldots \times S^1(a_N) \subset \C^N$, where~$a_i > 0$. Here,~$S^1(a) \subset \C$ denotes the circle bounding a disk of area~$a$.
}
\end{example}

Toric fibres are orbits with trivial stabilizer of the~$T^n$-action. This means that the torus action gives a canonical identification~$T^n \cong T(x)$ and
\begin{equation}
	\label{eq:lambda}
	\Lambda 
	= \ker(\exp \colon \ft \rightarrow T^n)
	= \pi_1(T(x))
	= H_1(T(x);\Z).
\end{equation}
Let us now discuss what happens to toric fibres under toric reduction. In general, let~$p \colon Z \rightarrow X$ be the quotient map of a symplectic reduction. If~$L \subset X$ is Lagrangian, then~$p^{-1}(L)$ is Lagrangian as well and we call it the \emph{lift} of~$L$. Conversely, any Lagrangian contained in~$Z$ is automatically invariant under the group action and projects to a Lagrangian in the reduced space. Adopting our notation from \S\ref{ssec:toricred}, let~$X'$ be a quotient obtained from~$X$ by toric reduction and let~$\iota(\Delta') \subset \Delta$ be the inclusion of the corresponding moment polytopes. Furthermore, we denote the reduction map by~$p \colon Z \rightarrow X'$ and the toric fibres in~$X$ by~$T(\cdot)$ and those in~$X'$ by~$T'(\cdot)$. 

\begin{proposition}
\label{prop:toricfibresred}
In the above notation, we have the following correspondence of toric fibres in~$X$ and~$X'$, 
\begin{equation}
	p^{-1}(T'(x)) = T(\iota(x)) \subset X, \quad
	x \in \Int \Delta'.
\end{equation}
\end{proposition}

\proof
This follows directly from the definition of the moment map~$\mu'$ on the quotient~$X'$. 
\proofend

In later sections, we will heavily use the second relative homotopy/homology groups of toric fibres, which is why we will discuss them here. Recall from~\eqref{eq:polytope} that the vectors~$\xi_i \in \Lambda$ are defined as orthogonal vectors to the facets of~$\Delta$. We prove the following well-known fact using the Delzant construction together with Proposition~\ref{prop:toricfibresred}.

\begin{proposition}
\label{prop:pi2}
Let~$(X,T(x))$ be a pair of a toric symplectic manifold and a toric fibre. Then~$\pi_2(X,T(x)) \cong \Z^N$, where~$N$ is the number of facets of~$\Delta$. Furthermore, there is a canonical basis~$D_1,\ldots,D_N \in \pi_2(X,T(x))$ bounding the classes~$\pp D_i = \xi_i \in \Lambda = \pi_1(T(x))$.
\end{proposition}

\proof
By the Delzant construction and Proposition~\ref{prop:toricfibresred}, the toric fibre~$T(x)$ lifts to a product torus~$T(\ell(x)) \subset \C^N$ under the reduction map~$p \colon Z \rightarrow X$, where~$N$ is the number of facets of~$\Delta$. Let~$\tilde{D}_1,\ldots,\tilde{D}_N$ be the obvious basis of~$\pi_2(\C^N,T(\ell(x)))$. Note that these can be chosen to lie in~$Z \subset \C^N$ since the image of~$Z$ under the moment map~$\mu_0$ is equal to the image of the embedding~$\ell$ from~\eqref{eq:ell}. Furthermore, reduction maps induce isomorphisms of relative homotopy groups, see for example the proof of~\cite[Proposition 3.2]{Smi19}. This shows that~$\pi_2(\C^N,T(\ell(x)))$ and~$\pi_2(X,T(x))$ are isomorphic, and we denote the image of~$\tilde{D}_i$ under the isomorphism by~$D_i$. In order to compute the boundary operator~$\pp$, consider the commutative diagram
\begin{equation}
\label{eq:homotopycd}
	\begin{tikzcd}
		\pi_2(\C^N,T(\ell(x))) \arrow{r}{\pp'} \arrow{d}{p_*}
		& \pi_1(T(\ell(x))) \arrow{d}{p_*} \\
		\pi_2(X,T(x)) \arrow{r}{\pp}
		& \pi_1(T(x)).
	\end{tikzcd}
\end{equation}
The boundary operator~$\pp'$ is an isomorphism mapping~$\tilde{D}_i$ to the $i$-th standard basis vector~$e_i$ and therefore it suffices to understand~$p_*$ on the fundamental group. Recall from the discussion surrounding~\eqref{eq:sestori} that~$p$ is equivariant in the sense that~$p(t.z)=\Xi(t).p(z)$ for all~$t \in T^N$ and~$z \in \C^N$. In the special case of the Delzant construction, one can easily check that~$\Xi_*(e_i)=\xi_i$ and thus this proves the last claim. \proofend

The homotopy long exact sequence for the pair~$(X,T(x))$ gives a short exact sequence, 
\begin{equation}
\label{eq:seshomotopy}
0 \rightarrow \pi_2(X) \rightarrow \pi_2(X,T(x)) \rightarrow \pi_1(T(x)) \rightarrow 0.
\end{equation}
Indeed, the higher homotopy grups of the torus vanish and toric manifolds are simply-connected. In homology (with integer coefficients) we obtain the same short exact sequence, 
\begin{equation}
\label{eq:seshomology}
0 \rightarrow H_2(X) \rightarrow H_2(X,T(x)) \rightarrow H_1(T(x)) \rightarrow 0.
\end{equation}
Indeed, the maps~$H_*(T(x)) \rightarrow H_*(X)$ are zero, since there is a contractible subset~$\Omega \subset X$ such that~$T(x) \subset \Omega \subset X$. Take for example~$\Omega = \mu^{\inv}(\Int \Delta \cup U)$, where~$U$ is a small neighbourhood of a vertex of~$\Delta$. There are obvious identifications of the respective groups in~\eqref{eq:seshomotopy} and~\eqref{eq:seshomology} which commute with the maps of these short exact sequences and thus we use homology and homotopy groups interchangeably.

Note that this discussion yields a very effective way to read off~$\pi_2(X) = H_2(X)$ from the moment polytope of a toric manifold. It is the kernel of~$\pp$, i.e.\ the lattice of integral relations among the vectors~$\xi_1,\ldots,\xi_N$ orthogonal to the facets of~$\Delta$. This in turn has a nice geometric interpretation in terms of the singular fibration structure of the moment map~$\mu \colon X \rightarrow \Delta$. Indeed, when moving from the interior of the moment polytope to the interior of a facet~$F_i$, the circle~$S^1(\xi_i) \subset T^n$ collapses, where by~$S^1(\xi_i)$ we have denoted the circle generated by the orthogonal vector~$\xi \in \Lambda \subset \ft$ to the facet~$F_i$. The canonical basis~$D_1,\ldots,D_N \in \pi_2(X,T(x))$ corresponds to the disks coming from these circles collapsing, which explains~$\pp D_i = \xi_i$. Furthermore, let~$\sum_i a_i D_i$ be an integral combination of such disks with~$\sum_i a_i \xi_i = 0$. The latter condition means that the corresponding concatenation of curves representing the~$\xi_i $ bound in the fibre~$T(x)$. Thus they define a homotopy class in~$\pi_2(X)$, which illustrates~$\ker \pp = \pi_2(X)$.

\section{Symmetric probes}
\label{sec:symmprobes}

Symmetric probes were first defined in~\cite{AbrBorMcD14}, where they were used to a different end. Let~$(X,\omega,\mu)$ be a toric symplectic manifold with moment polytope~$\Delta$. 

\begin{definition}
\label{def:symmprobe}
A \emph{symmetric probe}~$\sigma \subset \Delta$ is a reduction-admissible line segment, see Definition~\ref{def:redadm}.
\end{definition}

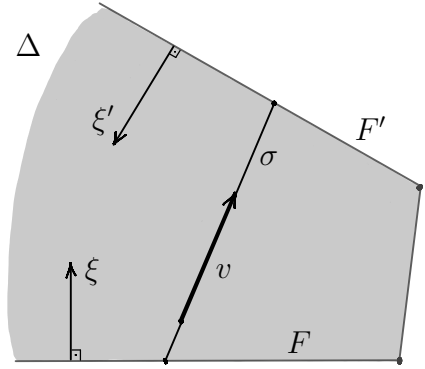
\begin{figure}
	\begin{center}
	\begin{tikzpicture}
		\node[inner sep=0pt] at (0,0)
    			{\includegraphics[trim={3.5cm 5cm 8cm 2.8cm},clip,scale=0.75]						{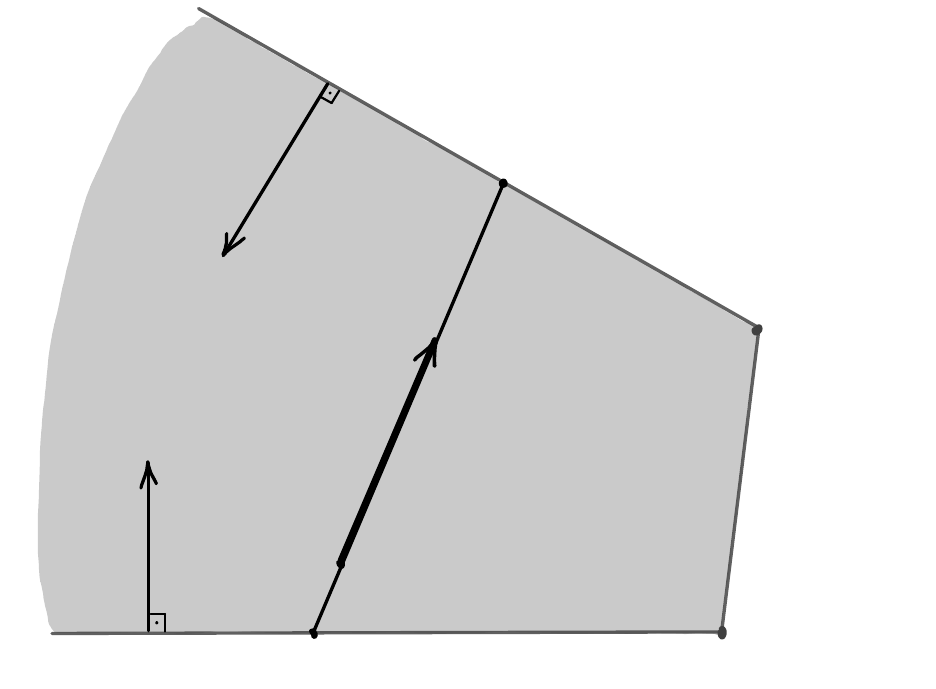}};
    		\node at (-3,2){$\Delta$};
    		\node at (-2,1){$\xi'$};
    		\node at (-2.15,-1){$\xi$};
    		\node at (1.55,0.9){$F'$};
    		\node at (0.6,-1.95){$F$};
    		\node at (-0.4,-1){$v$};
    		\node at (0.2,0.5){$\sigma$};
	\end{tikzpicture}
	\caption{A symmetric probe~$\sigma \subset \Delta$ and the surrounding notation.}
	\label{fig:4}
	\end{center}
\end{figure}

Let us unpack this definition and introduce some notation. By~$l \subset \ft^*$ we denote the line containing the symmetric probe~$\sigma$, by~$v \in \Lambda^*$ a primitive directional vector of~$l$ and by~$F$ and~$F'$ the facets of~$\Delta$ which~$\sigma$ intersects. We choose~$F,F'$ so that~$v$ points away from~$F$ and towards~$F'$. See Figure~\ref{fig:4} for an illustration of the set-up. Note that symmetric probes indeed do intersect facets, and not lower dimensional faces. Definition~\ref{def:symmprobe} implies that there is a basis of~$\Lambda^*$ contained in the unions~$l \cup F$ and~$l \cup F'$, respectively. This means that, locally, all intersections of symmetric probes with a facet are equivalent under integral affine transformations. After choosing a basis, we can work in~$(\R^n,\Z^n)$ and assume that
\begin{equation}
	\label{eq:inttransnf}
	v = e_n^*, \quad
	F = \Span_{\R}\{ e_1^*, \ldots , e_{n-1}^* \}.
\end{equation}
This follows from the fact that~$\GL(n;\Z)$ acts transitively on the set of bases of~$\Z^n$. We take~\eqref{eq:inttransnf} to be the normal form of an intersection of a symmetric probe with a facet. McDuff~\cite{McD11} calls these intersections \emph{integrally transverse} and we refer to her paper for a detailed discussion of this notion. In the above notation we have~$\langle v,\xi \rangle = -\langle v,\xi' \rangle = 1$, where~$\xi,\xi' \in \Lambda$ are the normal vectors to~$F$ and $F'$, respectively. By the normal form~\eqref{eq:inttransnf}, it follows that we can assume~$\xi = e_n$, which implies that~$\xi' = \sum_{i=1}^{n-1} k_ie_i - e_n $. The numbers~$k_1,\ldots,k_{n-1} \in \Z$ completely determine the toric structure of a neighbourhood of the symmetric probe~$\sigma$ and they are topological invariants of the torus bundle coming from the reduction map~$\mu^{-1}(\sigma) \rightarrow S^2$ appearing in the proof of Theorem~\ref{thm:symmprobe}.

\begin{theorem}
\label{thm:symmprobe}
Let~$\sigma \subset \Delta$ be a symmetric probe and~$x,y \in \sigma$ be a pair of points lying at equal distance to the boundary of~$\sigma$. Then the toric fibres~$T(x)$ and~$T(y)$ are Hamiltonian isotopic by a Hamiltonian isotopy inducing the map
\begin{equation}
	\label{eq:basicinvolution}
	\Phi_{\sigma} \colon H_1(T(x)) \rightarrow H_1(T(y)), \quad
	a \mapsto a + \langle v , a \rangle (\xi' - \xi)
\end{equation}
on the first homology of the toric fibres.
\end{theorem}

In particular, this proves Theorem~\ref{thm:mainA}. In~\eqref{eq:basicinvolution}, we have used the identification~$\Lambda = H_1(T(x)) = H_1(T(y))$ induced by the torus action. The map~$\Phi_{\sigma}$ is an involution and its~$(+1)$-Eigenspace is~$(n-1)$-dimensional and given by the complement~$\sigma^0=v^0 \subset \ft = \Lambda \otimes \R$. Its~$(-1)$-Eigenspace is spanned by~$\xi' - \xi$. Note also that~$\Phi_{\sigma}$ is uniquely determined by~$\sigma \subset \Delta$ and more precisely by an arbitrarily small neighbourhood of~$\sigma$ in~$\Delta$. Indeed, if we exchange~$\xi$ and~$\xi'$, then~$v$ changes its sign by our convention.

In case~$x=y$, we obtain an interesting corollary about the Hamiltonian monodromy group (see Definition~\ref{def:hmg}) of the corresponding toric fibre.

\begin{corollary}
\label{cor:symmprobemon}
Let~$x$ be the midpoint of a symmetric probe~$\sigma$. Then the Hamiltonian monodromy group~$\ch_{T(x)}$ contains the element~$\Phi_{\sigma}$. 
\end{corollary}

\proofof{Theorem~\ref{thm:symmprobe}}
Since~$\sigma \subset \Delta$ is reduction-admissible, we can perform toric reduction by Theorem~\ref{thm:toricred}. The reduced space is a copy of~$S^2$ with a standard symplectic form of total area equal to the integral affine length of~$\sigma$. Under the reduction, the fibres~$T(x), T(y)$ are mapped to a pair of circles~$S^1_x, S^1_y \subset S^2$ which are orbits of the residual Hamiltonian circle action on~$S^2$, see Proposition~\ref{prop:toricfibresred}. Since~$x,y$ are at equal distance to the boundary of~$\sigma$, the circles~$S_x^1,S^1_y$ bound disks of the same area and thus can be exchanged by a Hamiltonian isotopy~$\varphi$ on~$S^2$. Lift this Hamiltonian isotopy from~$S^2$ to~$X$ by lifting its Hamiltonian function by the reduction map to~$\mu^{-1}(\sigma)$ and extending it (for example by cut-off) to the total space. See for example~\cite[Lemma 3.1]{AbrMac13} or~\cite[Lemma 3.1]{Bre20} for details on lifting Hamiltonian isotopies.

Let us now compute the map induced by~$\varphi$ on~$\Lambda$. We work with homotopy groups here, but the problem is exactly the same in homology by the discussion in~\S\ref{ssec:toricfibres}. Let~$d_x,d_x' \in \pi_2(S^2,S^1_x)$ and~$d_y,d_y' \in \pi_2(S^2,S^1_y)$ be the generators of relative homotopy groups such that~$d_x,d_y$ contain the south pole and~$d_x',d_y'$ contain the north pole such that~$d_x + d_x' = d_y + d_y' = [S^2]$ for a chosen orientation on~$S^2$. The map~$\varphi_*$ induced by the Hamiltonian isotopy~$\varphi$ on relative homotopy groups satisfies~$\varphi_*d_x = d_y'$ and~$\varphi_*d_y = d_x'$. Furthermore, the map~$\Phi_{\sigma}$ is uniquely determined by the properties
\begin{equation}
	\label{eq:phisigma}
	\Phi_{\sigma}(\xi) = \xi', \quad
	\Phi_{\sigma}\vert_{v^0} = \id_{v^0},
\end{equation}
where~$v^0 \subset \Lambda$ denotes the elements on which~$v \in \Lambda^*$ vanishes. Indeed,~$\xi$ is transverse to~$v^0$ since~$\langle v, \xi \rangle = 1$. We show that the lift of~$\varphi$ satisfies~\eqref{eq:phisigma}, which proves the claim. The second property in~\eqref{eq:phisigma} follows from the~$K_{\sigma}$-equivariance of the lift of~$\varphi$ where~$K_{\sigma} = \exp \sigma^0$ is the complementary torus of the probe~$\sigma$. Indeed,~$K_{\sigma} \subset T^n$ is the subtorus with respect to which the symplectic reduction~$\mu^{-1}(\sigma) \rightarrow S^2$ is carried out, see also~\eqref{eq:comptorus} and the proof of Theorem~\ref{thm:toricred}, and thus any Hamiltonian isotopy lifted from the reduced space is equivariant with respect to this group action. For the first property in~\eqref{eq:phisigma}, note that the map~$p_* \colon \pi_2(\mu^{-1}(\sigma),T(x)) \rightarrow \pi_2(S^2,S^1_{x})$ induced by symplectic reduction is an isomorphism. See for example the proof of~\cite[Proposition 3.2]{Smi19}. Therefore~$\pi_2(\mu\inv(\sigma),T(x_{\sigma}))$ is generated by~$D_x,D'_x$ with~$\pi_*(D_x) = d_x$ and~$\pi_*(D'_x) = d'_x$ and similarly for~$T(y)$. This allows us to conclude that the lift of~$\varphi$ maps~$D_x$ to~$D_y'$ and~$D_y$ to~$D_x'$. Since~$\pp_xD_x = \pp_y D_y = \xi$ and~$\pp_y D_y = \pp_y D_y' = \xi'$ this finishes the proof.
\proofend

Note that we have actually computed the map induced on relative second homology, 
\begin{equation}
	H_2(X,T(x)) \rightarrow H_2(X,T(y)), \quad
	b \mapsto b + \langle v , \pp b \rangle (D' - D ),
\end{equation}
where~$D, D'$ denote the homology classes of the canonical basis in Proposition~\ref{prop:pi2} corresponding to~$F,F'$, respectively. Note also that the lift of~$\varphi$ in the proof of Theorem~\ref{thm:symmprobe} depends on the extension of the Hamiltonian function to~$X$ and is thus not uniquely defined by~$\varphi$. 

\begin{remark}
By choosing a suitable cut-off of the lifted Hamiltonian function in the proof of Theorem~\ref{thm:symmprobe}, one can choose the Hamiltonian isotopy to be supported in an arbitrarily small neighbourhood of~$\sigma \subset \Delta$.
\end{remark} 

\section{Chekanov invariants}
\label{sec:chekanovinv}

The main idea of this section is to use the Delzant construction to lift toric fibres of certain toric manifolds to product tori in some~$\C^N$ via Proposition~\ref{prop:toricfibresred} and to make use of the various results on product tori in~\cite{Che96}. In particular, this yields strong obstructions to the equivalence of toric fibres (Theorem~\ref{thm:mainB}) and their Hamiltonian monodromy (Theorem~\ref{thm:mainC}). As we shall discuss, similar results can be obtained \emph{by hand} (i.e.\ avoiding the lifting trick) via displacement energy and versal deformations, which comes in handy in case~$X$ cannot be seen as a toric reduction of~$\C^N$. However, we note that the approach by hand runs into the question of determining the displacement energy of toric fibres, which turns out to be very subtle in general, see for example the papers~\cite{AbrBorMcD14, McD11} for detailed discussions of the (qualitative) question of displaceability and~\cite[Section 3]{Bre20} for the quantitative question about displacement energy. In case~$X$ can be seen as a toric reduction of~$\C^N$, this question can be completely avoided by the lifting trick.

\begin{definition}
\label{def:redtype}
A toric symplectic manifold~$X$ is called \emph{of reduction type} if it can be obtained as a toric reduction of some~$\C^N$. 
\end{definition}

By the Delzant construction in~\S\ref{ssec:delzant}, all compact toric manifolds are of reduction type. The space~$X = \C \times S^2$, which will be discussed in~\S\ref{ssec:excs2}, is an example of a non-compact space which is of reduction type.

Before moving to the Chekanov invariants, let us point out the following.

\begin{remark}[Classification up to symplectomorphisms]
\label{rk:sympclass}
We focus on equivalence of toric fibres up to Hamiltonian diffeomorphisms. One may ask an analogue of Question~\ref{qu:class} for the group of symplectomorphisms. Note that a toric symplectic manifold~$X$ is simply connected whenever its moment polytope has at least one vertex, meaning that the distinction between the two classification questions is, at best, a question about connected components of~$\Symp(X,\omega)$. In fact, both classifications agree in all simply connected examples we consider in Section~\ref{sec:examples} of this paper. This is not always true, as the following example illustrates. Let~$X$ be the space obtained from~$\CP^2$ by three small toric blow-ups of the same size~$\varepsilon > 0$ at the vertices of the original moment triangle. The resulting symplectic manifold is toric and its moment polytope is a hexagon with three long and three short edges. Near each of the short edges, there is a non-displaceable toric fibre, as was proved in~\cite[\S 5.5]{Cho08}. In particular, these three non-displaceable fibres are not equivalent under Hamiltonian diffeomorphisms. However they are symplectomorphic. Indeed, their base points can be permuted by integral affine symmetries of the moment polytope and such symmetries lift to symplectomorphisms of the corresponding toric manifold, see for example~\cite[Lemma 4.3]{BreKimMoo19} for a proof of this well-known fact.

In particular, Conjecture~\ref{conj:main} is false for equivalence up to symplectomorphisms.
\end{remark}

\subsection{Equivalence of toric fibres}
\label{ssec:chekanovinv}

As we have seen in Example~\ref{ex:prodtori}, the product tori
\begin{equation}
	T(a) 
	= T(a_1,\ldots,a_N) 
	= S^1(a_1) \times \ldots \times S^1(a_N) \subset \C^N
\end{equation} 
are a special case of toric fibres. Chekanov has given a classification of product tori\footnote{Chekanov calls these tori \emph{elementary tori}.} up to symplectomorphism in~\cite[Theorem A]{Che96}. A complete set of invariants is given by
\begin{eqnarray}
d(a) &=& \min\{a_1,\ldots,a_N\}, \\
\#_d(a) &=& \#\{ i \in \{1,\ldots,N \} \,\vert\, a_i = d(a) \}, \\
\Gamma(a) &=& \Z\langle a_1 - d(a), \ldots, a_N - d(a) \rangle,
\end{eqnarray}
where we write~$a = (a_1,\ldots,a_N) \in \R_{>0}^N$. The first invariant is a positive real number and corresponds to the displacement energy\footnote{In the original paper, Chekanov uses the first Ekeland--Hofer capacity instead.}~$d(a) = e(\C^N,T(a))$. The second invariant is a positive integer less or equal to~$N$ (with equality if~$T(a)$ is monotone) which comes from versal deformations and displacement energy. As it turns out, versal deformations of product tori are given as the minimum of~$\#_d(a)$ linear functionals, and they contain no other information beyond this number. The third invariant~$\Gamma(a) \subset \R$ is a subgroup of~$\R$ generated by~$N - \#_d(a)$ elements and is a purely \emph{soft} invariant. In fact, it is the set of symplectic areas of disks with vanishing Maslov class~$m(\cdot)$. Note that in the case of~$\C^N$, the symplectic form has a primitive~$\lambda$ and thus we can express~$\Gamma(a)$ as
\begin{equation}
	\Gamma(a) = \left\{ \left. \int_{\gamma} \lambda \in \R \;\right\vert \; \gamma \in H_1(T(a)), \;  m(\gamma) = 0 \right\}.
\end{equation}

This invariant can be more explicitely expressed as~$\Gamma(a) = \Z\langle a_1 - d(a), \ldots , a_n - d(a) \rangle$.

\begin{theorem}[Chekanov]
\label{thm:chekanov}
The product tori~$T(a)$ and~$T(a')$ are symplectomorphic in~$\C^N$ if and only if
	\begin{equation}
		d(a) = d(a'), \quad
		\#_d(a) = \#_d(a'), \quad
		\Gamma(a) = \Gamma(a').
	\end{equation} 
\end{theorem}

Let us now get back to the case of toric fibres and prove Theorem~\ref{thm:mainB}. Recall from~\S\ref{ssec:introclass} that the \emph{Chekanov invariants} of a toric fibre~$T(x)\subset X$ are defined in terms of the integral affine distances~$\ell(x) = (\ell_1(x),\ldots,\ell_N(x))$ of the the point~$x$ to the facets of~$\Delta$,
\begin{eqnarray}
\label{eq:ell}
d(x) &=& \min\{\ell_1(x),\ldots,\ell_N(x)\}, \\
\#_d(x) &=& \#\{ i \in \{1,\ldots,N \} \,\vert\, \ell_i(x) = d(x) \}, \\
\Gamma(x) &=& \Z\langle \ell_1(x) - d(x), \ldots, \ell_N(x) - d(x) \rangle.
\end{eqnarray}

\proofof{Theorem~\ref{thm:mainB}}
We prove the result for all toric manifolds of reduction type, see Definition~\ref{def:redtype}. Let~$X$ be a toric manifold of reduction type and~$T(x),T(x') \subset X$ be toric fibres which are equivalent under Hamiltonian isotopies. Recall from~\S\ref{ssec:delzant} that we may view~$X$ as a toric reduction of~$\C^N$, where the inclusion map of the moment polytope~$\Delta$ of~$X$ into~$\R^N$ is given by the map~$\ell(x) = (\ell_1(x),\ldots,\ell_N(x))$. Furthermore, the toric fibre~$T(x)$ lifts to the product torus~$T(\ell(x))$ in~$\C^N$ by Proposition~\ref{prop:toricfibresred} and similarly for~$T(x')$. Since~$T(x)\cong T(x')$, we obtain that~$T(\ell(x))\cong T(\ell(x'))$. Indeed, Hamiltonian isotopies can be lifted through symplectic reductions by lifting the corresponding Hamiltonian function and extending it to~$\C^N$ by cut-off. It is easy to see that any such lift will map the lift of~$T(x)$ to the lift of~$T(x')$. Theorem~\ref{thm:mainB} now follows from Theorem~\ref{thm:chekanov}. 
\proofend

The Chekanov invariants are not complete, as the following example illustrates. 

\begin{example}
\label{ex:cp2notcomplete}
Let~$\CP^2$ the complex projective plane equipped with the toric structure described in Example~\ref{ex:cp2} and with moment polytope~$\Delta$, and set
\begin{equation}
	x = \left( -\frac{5}{10} , -\frac{2}{10} \right), \quad
	x' = \left( -\frac{5}{10} , \frac{1}{10} \right) \in \Delta.
\end{equation}
Since $\ell(x)=(1+x_1,1+x_2,1-x_1-x_2)$, we obtain
\begin{equation}
	\ell(x) = \left( \frac{5}{10} , \frac{8}{10}, \frac{17}{10} \right), \quad
	\ell(x') = \left( \frac{5}{10} , \frac{11}{10}, \frac{14}{10} \right) \in \R^3_{\geqslant 0}.
\end{equation}
By the classification of toric fibres in~$\CP^2$ from~\cite[Proposition 7.1]{SheTonVia19}, see also~\S\ref{ssec:cp2}, the fibres~$T(x)$ and~$T(x')$ are not Hamiltonian isotopic. However, their Chekanov invariants agree. Indeed, we find
\begin{equation}
d(x)=d(x')=\frac{1}{2}, \quad
\#_d(x)=\#_d(x')=1, \quad
\Gamma(x)=\Gamma(x')=\Z\left\langle \frac{3}{10} \right\rangle.
\end{equation}
\end{example}

\subsection{Hamiltonian monodromy}
\label{ssec:hammon}

Let~$\phi \in \Ham(X,\omega)$ be a Hamiltonian diffeomorphism of a toric manifold~$(X,\omega)$ mapping a toric fibre~$T(x)$ to a toric fibre~$T(x')$. Then one can consider the map induced on relative second homology, 
\begin{equation}
	\phi_* \colon H_2(X,T(x)) \rightarrow H_2(X,T(x')).
\end{equation}
We call this map \emph{ambient monodromy}. In the same vein as in~\S\ref{ssec:chekanovinv}, we derive obstructions to which maps~$\phi_*$ can be obtained in this way by using the Delzant construction to lift Hamiltonian isotopies. Note that by setting~$x = x'$ and by projecting to the first homology (see~\eqref{eq:seshomology}), we can extract information about the Hamiltonian monodromy question as a special case. 

The key result by Chekanov is~\cite[Theorem 4.5]{Che96}.

\begin{theorem}[Chekanov]
\label{thm:chekanov2}
Let~$T(a),T(a') \subset \C^N$ be product tori. An isomorphism
\begin{equation}´
\Phi \colon H_1(T(a)) \rightarrow H_1(T(a'))
\end{equation}
can be realized as~$(\phi\vert_{T(a)})_* = \Phi$ by a symplectomorphism~$\phi \in \Symp(X,\omega)$ mapping~$T(a)$ to~$T(a')$ if and only if the following conditions hold
\begin{equation}
	\Phi(\mathcal{D}(a)) = \mathcal{D}(a'), \quad
	\Phi^* m_{T(a')} = m_{T(a)}, \quad
	\Phi^* \sigma_{T(a')} = \sigma_{T(a)}.
\end{equation}
\end{theorem}
Here~$m_{T(a)} \in H^1(T(a);\Z)$ and~$\sigma_{T(a)} \in H^1(T(a);\R)$ are the Maslov class and the symplectic area class, respectively. By~$\cd(a) \subset H_1(T(a))$ we denote the set of \emph{distinguished classes}. In the standard basis~$e_1,\ldots,e_N \in H_1(T(a))$ the basis vector~$e_i$ is called a \emph{distinguished class} if the corresponding component in~$a = (a_1,\ldots,a_N)$ is minimal, i.e.\ if ~$a_i = d(a)$.

Let us now move to toric fibres. Recall from Proposition~\ref{prop:toricfibresred} that for any toric fibre~$T(x) \subset X$, the relative second homology~$H_2(X,T(x))$ has a canonical basis~$D_1,\ldots,D_N$, where~$D_i$ corresponds to the~$i$-th facet of the moment polytope~$\Delta$ of~$X$. 

\begin{definition}
Let~$T(x)\subset X$ be a toric fibre. The \emph{distinguished classes} of~$T(x)$ are the elements of the set
	\begin{equation}
		\cd(x) = \{D_i \, \vert \, \ell_i(x) = d(x)\} \subset H_2(X,T(x)), 
	\end{equation}
i.e.\ elements of the canonical basis for which the distance of~$x \in \Int \Delta$ to the corresponding facet of~$\Delta$ is minimal.
\end{definition}

Recall that there is a canonical inclusion~$H_2(X) \subset H_2(X,T(x))$, meaning that there is a distiguished subspace which is independent of the choice of~$x$. We prove the following.

\begin{theorem}
\label{thm:ambientmonodromy}
Let~$T(x),T(x') \subset X$ be toric fibres in a compact toric manifold~$X$ such that there exists a Hamiltonian diffeomorphism~$\phi \in \Ham(X,\omega)$ mapping~$T(x)$ to~$T(x')$. Then the induced map 
\begin{equation}
	\phi_* \colon H_2(X,T(x)) \rightarrow H_2(X,T(x')) 
\end{equation} 
on relative homology groups satisfies
\begin{equation}
	\label{eq:exmonodromyconstraints}
	\phi_*(\cd(x)) = \cd(x'), \quad
	\phi^*m_{T(x')} = m_{T(x)}, \quad
	\phi^*\sigma_{T(x')} = \sigma_{T(x)}, \quad
	\phi_*\vert_{H_2(X)} = \id.
\end{equation}
\end{theorem}

\proof
The second and third identity in~\eqref{eq:exmonodromyconstraints} are general facts about the Maslov and the symplectic area class. The last identity is straightforward since Hamiltonian diffeomorphisms are isotopic to the identity on~$X$ and hence the map~$\phi_* \colon H_*(X) \rightarrow H_*(X)$ is the identity.  For the first identity in~\eqref{eq:exmonodromyconstraints}, we again use the Delzant construction together with lifting the Hamiltonian isotopy. The following groups are canonically isomorphic,
\begin{equation}
	\label{eq:canisom}
	H_2(X,T(x)) 
	\cong H_2(\C^N,T(\ell(x)))
	\cong H_1(T(\ell(x))),
\end{equation}
see the proof of Proposition~\ref{prop:pi2}, where this is proved for the corresponding (relative) homotopy groups. Thus the map~$H_1(T(\ell(x))) \rightarrow H_1(T(\ell(x')))$ induced by the lifted Hamiltonian diffeomorphism is conjugate to~$\phi_*$ by the canonical isomorphism~\eqref{eq:canisom}. It is easy to see that the distinguished classes~$\cd(x)$ of~$T(x)$ are by definition mapped under~\eqref{eq:canisom} to the distinguished classes~$\cd(\ell(x))$ of the product torus~$T(\ell(x))$ and thus the first identity in~\eqref{eq:exmonodromyconstraints} follows from Theorem~\ref{thm:chekanov2}. 
\proofend

It seems reasonable to guess that these constraints are sufficient. More precisely, note that there is a canonical identification~$H_2(X,T(x)) = H_2(X,T(x'))$ for any two points~$x,x' \in \Int \Delta$. Then we conjecture the following.

\begin{conjecture}
An isomorphism~$\Phi \in \Aut H_2(X,T(x))$ can be realized as ambient monodromy of a Hamiltonian diffeomorphism mapping~$T(x)$ to~$T(x')$ if and only if the identities in~\eqref{eq:exmonodromyconstraints} hold.
\end{conjecture}

We show that this conjecture holds in all examples discussed in Section~\ref{sec:examples}. In fact, we use the ambient monodromy and Theorem~\ref{thm:ambientmonodromy} to classify toric fibres and determine the Hamiltonian monodromy groups in these examples. The area class~$\sigma_{T(x)}$ determines~$x$ and hence proving this conjecture gives, in particular, an answer to Question~\ref{qu:class}.

Let us now move to the ordinary Hamiltonian monodromy group of toric fibres, see Definition~\ref{def:hmg}. To derive information about~$\ch_{T(x)}$ from Theorem~\ref{thm:ambientmonodromy}, fix~$x = x'$ and let~$\phi \in \Ham(X,\omega)$ be a Hamiltonian isotopy such that~$\phi(T(x)) = T(x)$. Note that the ambient monodromy~$\phi_*$ determines the map~$(\phi\vert_{T(x)})_* \in \Aut H_1(T(x))$ by the short exact sequence~\eqref{eq:seshomology}.\medskip

\proofof{Theorem~\ref{thm:mainC}}
Any element in the Hamiltonian monodromy group~$\ch_{T(x)}$ comes from an ambient monodromy element~$\phi_*$ by~\eqref{eq:seshomology} and hence the theorem follows directly from Theorem~\ref{thm:ambientmonodromy} where the set of distinguished classes in~$H_1(T(x))$ is given by
\begin{equation}
	\label{eq:disth1}
	\pp \cd(x)
	= \{\xi_i \, \vert \, \ell_i(x) = d(x)\}
	\subset H_1(T(x)),
\end{equation}
where~$\xi_i \in \ft \cong H_1(T(x))$ is a primitive defining vector of the~$i$-th facet of~$\Delta$. Indeed, recall from Proposition~\ref{prop:pi2} that the boundary of a canonical basis element~$D_i$ is~$\xi_i$.
\proofend

It follows from Theorem~\ref{thm:mainC} that if the distinguished classes span the lattice~$H_1(T(x))$, then~$\ch_{T(x)}$ is a subgroup of the group of permutations on~$\#_d(x)$ elements. In particular, the Hamiltonian monodromy group is finite in this case. See also~\cite[Theorem 1]{AugSmiWor22}. In contrast, we shall see that the Hamiltonian monodromy group is infinite in some examples, see~\S\ref{ssec:excs2}, \S\ref{ssec:c2ts1} and~\S\ref{ssec:ts1s2}. The number~$\#_d(x)$ is maximal in case~$T(x)$ is the monotone toric fibre of a (monotone) toric manifold~$X$. In that case, we obtain the obstructive statement of~\cite[Theorem 2]{AugSmiWor22} for the group of Hamiltonian diffeomorphisms as a special case of Theorem~\ref{thm:ambientmonodromy}. 

\begin{corollary}
\label{cor:smith}
Let~$T(x) \subset X$ be a monotone toric fibre. Then any element in~$\ch_{T(x)}$ acts as a permutation on the set~$\{\xi_1, \ldots, \xi_N\}$ of defining vectors of the polytope and the corresponding ambient monodromy acts as the identity on~$H_2(X)$.
\end{corollary}

\subsection{Displacement energy and versal deformations of toric fibres}
\label{ssec:displacementen}

In this subsection, we discuss obstructions for the equivalence of toric fibres and their Hamiltonian monodromy relying on versal deformations instead of the lifting trick employed in the proofs of Theorems~\ref{thm:mainB} and \ref{thm:ambientmonodromy}. This comes in handy in cases where~$X$ cannot be seen as a toric reduction of some~$\C^N$, and we will use them in~\S\ref{ssec:c2ts1} and~\S\ref{ssec:ts1s2}. Note that the direct approach by versal deformations has the drawback that it requires a computation of the displacement energy of toric fibres, at least on an open dense subset. See Assumption~\ref{ass:detoric}. 

Let us briefly discuss displacement energy and versal deformations. We refer to~\cite{Che96} and especially~\cite{CheSch10, CheSch16} for more details. The displacement energy of a compact subset~$A \subset (X,\omega)$ is defined as the infimum of the Hofer norm taken over all Hamiltonian isotopies displacing~$A$ from itself, 
\begin{equation}
	e(X,A) 
	=
	\inf\{\Vert H \Vert \, \vert \, \phi_1^H(A) \cap A = \varnothing \},
\end{equation}
and by convention~$e(X,A) = \infty$ if the infimum is taken over the empty set. The displacement energy is a symplectic invariant and we will use it only in case~$A$ is a Lagrangian. 

For a compact Lagrangian~$L \subset X$, Chekanov introduced a way to strengthen a given symplectic invariant by looking at the invariant on Lagrangian neighbours of~$L$. This is called \emph{versal deformation} of~$L$. Perturbing~$L$ in a Weinstein neighbourhood, we find that nearby Lagrangians correspond to graphs of closed one-forms on~$L$. Furthermore, we can associate to every such perturbation an element in~$H^1(L;\R)$, by taking its (Lagrangian) flux. Two such perturbations are Hamiltonian isotopic (with support in the Weinstein neighbourhood of~$L$) if and only if they map to the same element in~$H^1(L;\R)$. Thus we obtain a continuous bijection between locally supported Hamiltonian isotopy classes of Lagrangian neighbours of~$L$ and a neighbourhood of the origin of~$H^1(L;\R)$. As the flux description suggests, this correspondence is independent of the chosen Weinstein neighbourhood.

We may post-compose any symplectic invariant with the map from~$U \subset H^1(L;\R)$ to classes of nearby Lagrangians. Here, we use displacement energy to obtain a function~$U \rightarrow \R \cup \{\infty\}$. By taking its germ, we obtain,
\begin{equation}
	\label{eq:degerm}
	\ce_L \colon H^1(L;\R) \rightarrow \R \cup \{ \infty \}.
\end{equation}

\begin{definition}
We call the function~\eqref{eq:degerm} the \emph{displacement energy germ} of~$L \subset X$.
\end{definition}

The displacement energy germ is a symplectic invariant in the sense that if~$\phi \in \Symp(X,\omega)$, then 
\begin{equation}
\label{eq:vdinv}
	\ce_L \circ \phi\vert_L^* = \ce_{\phi(L)},
\end{equation}
where~$\phi\vert_L^*$ is the transpose of the isomorphism~$(\phi\vert_L)_* \colon H_1(L) \rightarrow H_1(\phi(L))$. In particular, this can be used to derive obstructions to Hamiltonian monodromy. 

\begin{proposition}
\label{prop:vdmonodormy}
Let~$L \subset X$ be a compact Lagrangian submanifold. If~$\Phi \in \ch_L$ is an element in the Hamiltonian monodromy group, then~$\ce_L \circ \Phi^* = \ce_L$. 
\end{proposition}

Let us discuss this in more detail in the special case where~$L= T(x) \subset X$ is a toric fibre of a toric manifold~$(X,\omega)$. Coming up with a versal deformation of toric fibres is straightforward. Indeed, a versal deformation of~$T(x)$ is obtained by varying the base point,~$a \mapsto T(x+a)$ for small enough~$a$, where we identify~$H^1(T(x);\R) \cong \ft^*$ as usual via the~$T^n$-action. Thus the crucial point in computing~$\ce_{T(x)}$ is finding the displacement energy of toric fibres~$e(X,T(x))$ as a function of~$x \in \Int(T(x))$. Let us make the following assumption. 

\begin{assumption}
\label{ass:detoric}
On an open and dense subset of the moment polytope~$\Delta$, we assume that 
	\begin{equation}
	\label{eq:detoric}
		e(X,T(x)) 
		= d(x)
		= \min \{ \ell_1(x), \ldots, \ell_N(x) \}.
	\end{equation}
\end{assumption}
Here,~$d(\cdot)$ denotes the integral affine distance to the boundary of~$\Delta$ as in~\eqref{eq:ell}. Recall that the functionals~$\ell_i(\cdot) = \langle \cdot , \xi_i \rangle + \lambda_i$ measure the integral affine distance of~$x$ to the~$i$-th facet of~$\Delta$. Let~$f,g$ be two functions defined on a vector space~$V$. Since equalities on open and dense subsets will come up quite often and are in fact sufficient for our purposes, we write~$f \simeq g$ if~$f$ and~$g$ agree on an open and dense subset of~$V$. 

Let us briefly discuss why Assumption~\ref{ass:detoric} is reasonable. First, we note that the inequality~$e(X,T(x)) \geqslant d(x)$ holds whenever~$X$ is compact toric, and more generally, whenever~$X$ can be seen as the toric reduction of some~$\C^N$. This follows again from toric reduction and the lifting trick, see also~\cite[\S 3.2]{Bre20}. Indeed, if~$T(x) \subset X$ can be displaced with energy~$e$, then so can the corresponding product torus~$T(\ell(x)) \subset \C^N$ obtained by Proposition~\ref{prop:toricfibresred}. The displacement energy of the latter is precisely given by~$d(x)= \min \{ \ell_1(x), \ldots, \ell_N(x) \}$. Although this inequality may fail to be sharp (for example for non-dispaceable tori), in all the examples we know of, it fails only on the complement of an open dense subset, meaning that Assumption~\ref{ass:detoric} still holds. Furthermore, the assumption holds for all compact \emph{monotone} toric symplectic manifolds of dimension~$\leqslant 18$ as was checked computationally. The monotone case in arbitrary dimension is related to the so-called Ewald conjecture. See~\cite{McD11} or~\cite[\S 3.4]{Bre20} for a detailed discussion. The following proposition is~\cite[Proposition 4.3]{Bre20}.

\begin{proposition}
\label{prop:vdtoricfibres}
Under Assumption~\ref{ass:detoric}, the displacement energy germ of~$T(x)$ is given by 
\begin{equation}
	\ce_{T(x)}(a) 
	\simeq \min_{i \in I(x)} \{ \ell_i(x+a) \},
\end{equation}
where~$I(x) \subset \{1,\ldots,N\}$ is the subset of indices for which~$\ell_i(x)$ is minimal. 
\end{proposition}

Under Assumption~\ref{ass:detoric}, we can prove the symplectically hard part of Theorem~\ref{thm:mainB} and a weaker form of the hard part of Theorem~\ref{thm:ambientmonodromy}, where ambient monodromy is replaced by the map induced on first homology. 

\begin{theorem}
\label{thm:alternative}
Let~$X$ be a toric manifold for which Assumption~\ref{ass:detoric} holds. Let~$\phi$ be a Hamiltonian diffeomorphism mapping a toric fibre~$T(x)$ to a toric fibre~$T(x')$. Then we have
\begin{equation}
	\label{eq:alt}
	d(x) = d(x'), \quad
	\#_d(x) = \#_d(x').
\end{equation}
Furthermore, the map~$(\phi\vert_{T(x)})_* \colon H_1(T(x)) \rightarrow H_1(T(x'))$ acts by a permutation on distinguished classes,~$(\phi\vert_{T(x)})_*\cd(x) = \cd(x')$.
\end{theorem}

\proof
Let~$U \subset \Int \Delta$ be an open dense subset such that~\eqref{eq:detoric} holds for all~$x \in U$. For~$x,x' \in U$, we have~$d(x) = d(x')$. If~$x \notin U$ or~$x' \notin U$, use Proposition~\ref{prop:vdtoricfibres} to see that~$\min_{i \in I(x)} \{ \ell_i(x+a) \} \simeq \min_{i \in I(x')} \{ \ell_i(x'+a) \}$. Thus these two continuous functions of~$a$ are actually equal near~$a = 0$, and they yield~$d(x)$ and~$d(x')$, respectively, when evaluated at~$a =0$. The second invariance property in~\eqref{eq:alt} similarly follows from~\eqref{eq:vdinv} and Proposition~\ref{prop:vdtoricfibres} by noting that~$\#_d(x) = \#I(x)$. The claim about~$(\phi\vert_{T(x)})_*$ follows from~\eqref{eq:vdinv} and Proposition~\ref{prop:vdtoricfibres}. Indeed, recall that the distinguished classes of~$T(x)$ are the vectors~$\xi_i$ for which the corresponding~$\ell_i$ is minimal, see~\eqref{eq:disth1}.
\proofend

To illustrate that the methods of this paragraph can be applied to a broader set of examples than toric fibres, we include the following example.

\begin{example}[Vianna tori in~$\CP^2$]
\label{ex:vianna}
Using Proposition~\ref{prop:vdmonodormy}, one can show that all Vianna tori in~$\CP^2$, except for the first and the second one, have trivial Hamiltonian monodromy groups. The Vianna tori in~$\CP^2$ form a countable family of monotone Lagrangian tori which are not pairwise symplectomorphic. They are in bijection with so-called \emph{Markov triples}, i.e.\ triples of natural numbers solving the \emph{Markov equation}. We refer to~\cite{Via16} for a detailed description. We denote the Vianna torus corresponding to a Markov triple~$(a,b,c)$ by~$T(a,b,c) \subset \CP^2$. This torus appears as a monotone fibre of an \emph{almost toric fibration} of~$\CP^2$ with base diagram given by a certain triangle~$\Delta_{a,b,c}$. On an open and dense subset of a neighbourhood of the origin in~$H^1(T(a,b,c);\R)$, the displacement energy germ~$\ce_{T(a,b,c)}$ has level sets given by scalings of~$\pp\Delta_{a,b,c}$. This means that the versal deformation \emph{sees} the corresponding almost toric base diagram, and thus the integral affine equivalence class of~$\Delta_{a,b,c}$ is an invariant of~$T(a,b,c)$. In particular, this can be used to distiguish the Vianna tori as was noted by Chekanov--Schlenk in private communications. For a proof of this claim, see the forthcoming paper~\cite{BreMikSch21}. 

Using Proposition~\ref{prop:vdmonodormy}, we note that a necessary condition for~$T(a,b,c)$ to admit non-trivial monodromy is that the corresponding almost toric base diagram admits some integral affine symmetry. Such a symmetry can only exist if at least two vertices are of the same integral affine type, i.e.\ if the same Markov number appears at least twice in the same triple. This is only the case for~$(1,1,1)$ and~$(1,1,2)$. The former is the Clifford torus which has Hamiltonian monodromy group isomorphic to the dihedral group~$D_6$ and the latter is the first non-trivial Vianna torus~$T(1,1,2)$ having monodromy group isomorphic to~$\Z_2$. For all other Vianna tori, we obtain~$\ch_{T(a,b,c)} = \{1\}$. In particular, the Hamiltonian monodromy group does not contain enough information to distinguish Vianna tori.
\end{example}

\section{Examples}
\label{sec:examples}

In Subsections~\S\ref{ssec:s2s2}--\ref{ssec:ts1s2}, we classify toric fibres and determine their Hamiltonian monodromy in some examples. With the exception of~$\C^2 \times T^*S^1$, our examples are four-dimensional. This comes from the fact that the classification question in dimensions~$\geqslant 6$ is qualitatively very different -- provided the moment polytope has at least one vertex. Indeed, in that case, there are toric fibres~$T(x)$ for which~$\fH_x$ has accumulation points, see Corollary~\ref{cor:accumulation}. 

The proofs of the results of this section all follow the same pattern. Equivalences and monodromy elements are constructed by symmetric probes. The main ingredients for the obstructive side are Theorems~\ref{thm:mainB} and~\ref{thm:ambientmonodromy} applied to the ambient monodromy map~$\phi_* \colon H_2(X,T(x)) \rightarrow H_2(X,T(x'))$ induced by a Hamiltonian diffeomorphism~$\phi$. The conceptual reason why constraints on ambient monodromy give constraints on equivalences of toric fibres is the observation that the symplectic area class of a toric fibre determines~$x \in \Delta$. These methods probably apply to most four-dimensional toric manifolds, with the computational complexity increasing with the number of edges of the moment polytope. The examples we chose are diverse in the sense that~$S^2 \times S^2$ and~$\CP^2$ are compact toric and thus the Delzant construction can be used directly;~$\C \times S^2$ is non-compact, but still a toric reduction of~$\C^3$; the spaces~$\C^2 \times T^*S^1$ and~$T^*S^1 \times S^2$ are non-compact and cannot be seen as toric reductions of any~$\C^N$. However, the latter is a toric reduction of the former. In the case of the former, we apply the direct methods from~\S\ref{ssec:displacementen}. Note also that the spaces~$\C^2 \times T^*S^1$ and~$T^*S^1 \times S^2$ are not simply-connected, whence the classification up to symplectomorphisms is drastically different from the classification up to Hamiltonian diffeomorphisms. 

In~\S\ref{ssec:chekrev}, we revisit Chekanov's classification result and prove Conjecture~\ref{conj:main} for~$\C^n$. In~\S\ref{ssec:arbitrary}, we collect some remarks on how to construct symmetric probes in arbitrary toric manifolds.

Let us point out that all monodromy results for \emph{monotone} toric fibres in this section also follow from the methods developed in~\cite{AugSmiWor22}.

\subsection{The case of monotone $X = S^2 \times S^2$}
\label{ssec:s2s2}

Let~$S^2 \times S^2$ be equipped with the monotone product symplectic structure~$\omega = \omega_{S^2} \oplus \omega_{S^2}$, where~$\omega_{S^2}$ is the area form with normalization~$\int_{S^2} \omega_{S^2} = 2$. Then the corresponding moment polytope is given by the square~$\Delta = [-1,1] \times [-1,1]$. There are probes with four different directional vectors. The probes with~$v = e_1^*, e_2^*$ are admissible everywhere in the interior of the polytope. The probes with~$v = e_1^* + e_2^*$ and~$v = e_1^* - e_2^*$ are admissible everwhere except for the two main diagonals of the square. 
Note that the equivalences of toric fibres generated by these probes can also be read off from the symmetries of~$\Delta = [-1,1] \times [-1,1]$. Let us turn to the classification of toric fibres. 

\begin{proposition}
\label{prop:s2s2class}
The classification of toric fibres of monotone~$S^2 \times S^2$ is given by
	\begin{equation}
		\fH_x = \{ (\pm x_1 , \pm x_2) , (\pm x_2, \pm x_1) \}, \quad x = (x_1,x_2) \in \Int \Delta.
	\end{equation}
\end{proposition}

\begin{figure}
	\begin{center}
	\begin{tikzpicture}
		\node[inner sep=0pt] at (0,0)
    			{\includegraphics[trim={5cm 6.5cm 7cm 0.8cm},clip,scale=0.8]						{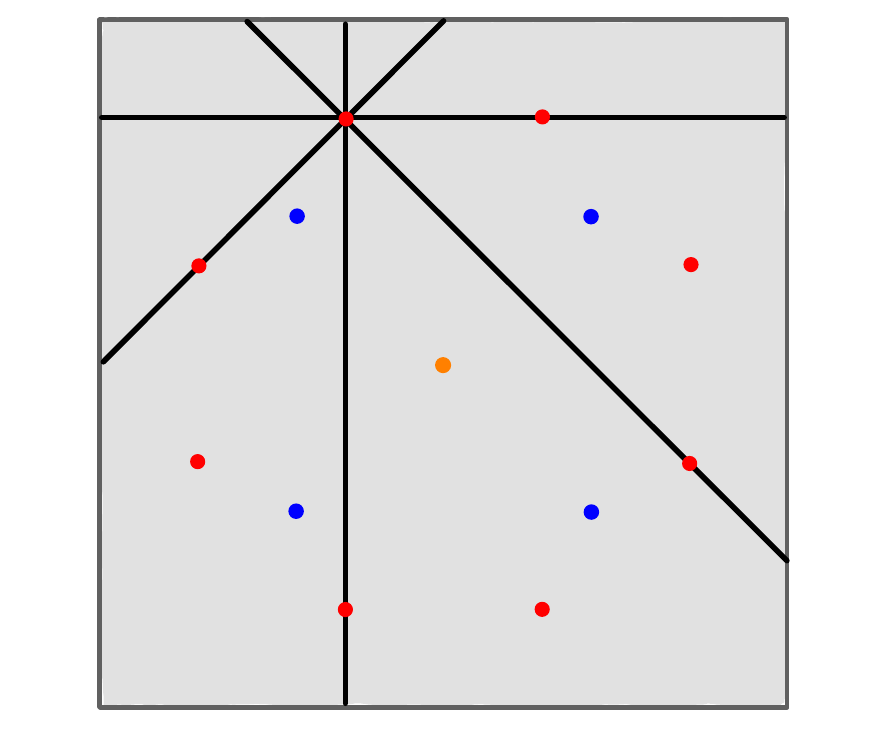}};
	\end{tikzpicture}
	\caption{Some symmetric probes in the monotone~$S^2\times S^2$. Points of the same colour denote equivalent fibres.}
	\label{fig:6}
	\end{center}
\end{figure}

Note that the sets~$\fH_x$ contain eight elements if~$x_1 \neq x_2$ and both are non-zero, four elements if~$x_1 = x_2$ or if one of the~$x_i$ is zero, and one element (the monotone fibre) if~$x_1 = x_2 = 0$. See Figure~\ref{fig:6}. 
\smallskip

\proofof{Proposition~\ref{prop:s2s2class}}
The constructive side follows either from the symmetric probes listed above or from the symmetries of~$\Delta$. For the obstructions, we will use Chekanov's invariants as expressed in Theorems~\ref{thm:mainB} and~\ref{thm:ambientmonodromy}. Let~$T(x)=T(x_1,x_2)$. By Chekanov's first invariant from Theorem~\ref{thm:mainB}, we can restrict our attention to the set~$\fD_x$ lying at distance~$d(x)$ to the boundary of~$\Delta$. This set is the boundary of a square of size~$2(1-d(x))$ and it is stratified by the second Chekanov invariant. Indeed, we have~$\#_d(x) = 2$ whenever~$x$ is a vertex of~$\fD_x$ and~$\#_d(x) = 1$ elsewhere on~$\fD_x$. This means that we are left with proving the result on the interior of the four edges of~$\fD_x$. Using the symmetries of~$\fD_x \subset \Delta$, we can restrict our attention to the segment~$[0,x_2) \times \{x_2\}$ since this is a fundamental domain for~$\fD_x$ under the symmetries. 

\textbf{Claim:}
If~$T(x_1,x_2) \cong T(x_1',x_2)$ for some~$x_1, x'_1 \in [0,x_2)$ then~$x_1 = x_1'$.

\ni
We use Theorem~\ref{thm:ambientmonodromy} to prove the claim. Suppose there is~$\phi \in \Ham(X,\omega)$ mapping~$T(x_1,x_2)$ to~$T(x_1',x_2)$. This induces
\begin{equation}
	\phi_* \colon H_2(X,T(x_1,x_2)) \rightarrow H_2(X,T(x_1',x_2)).
\end{equation}

Let~$D_1,D_2,D_3,D_4$ be the canonical basis of~$H_2(X,T(x_1,x_2))$ as in Proposition~\ref{prop:pi2} where~$D_1$ is the disk corresponding to the facet~$\{1\} \times [-1,1]$ and the remaining ones are ordered in the anti-clockwise direction. Let~$D_1',D_2',D_3',D_4'$ be the corresponding basis elements for~$H_2(X,T(x_1',x_2))$. The distinguished classes are~$\cd(x_1,x_2) =\{D_2\}$ and~$\cd(x_1',x_2) = \{D_2'\}$, meaning that Theorem~\ref{thm:ambientmonodromy} yields~$\phi_* D_2 = D_2'$. Set
\begin{equation}
	\phi_* D_1 = a_1D_1' + a_2D_2' + a_3D_3' + a_4D_4', \quad
	a_i \in \Z.
\end{equation}
Since the Maslov class is preserved, we obtain~$a_1 + a_2 + a_3 + a_4 = 1$, meaning that~$a_4 = 1 - a_1 - a_2 - a_3$. 
Since the induced map~$(\phi\vert_{T(x_1,x_2)})_*$ is invertible, we deduce that~$\det(\pp\phi_*D_1,\pp\phi_*D_2) = \pm 1$ which yields~$a_1 - a_3 = \pm 1$. Preservation of symplectic area,~$\int_{D_1} \omega = \int_{\phi_*D_1} \omega$ yields
\begin{equation}
	1 - x_1 = a_1(1-x_1') + a_2(1-x_2) + a_3(1+x_1') + a_4(1+x_2),
\end{equation}
since the areas of the disks~$D_i'$ are just given by the distances of~$(x_1',x_2)$ to the respective facets of~$\Delta$.
In case~$a_1 - a_3 = +1$, we use the above relations on the~$a_i$ to find~$a_3 = a_1 - 1$ and~$a_4 = 2 - 2a_1 - a_2$ and thus,
\begin{equation}
	\label{eq:xdiff1}
	x_1 - x'_1 = 2x_2(a_1 + a_2 - 1) \in 2x_2 \Z. 
\end{equation}
Since~$\vert x_1 - x'_1 \vert < x_2$, we conclude~$x_1 = x_1'$. In case~$a_1 - a_3 = -1$, we find by the same reasoning,
\begin{equation}
	\label{eq:xdiff2}
	x_1 + x_1' = 2x_2(a_1 + a_2) \in 2x_2 \Z.
\end{equation}
Since~$0 \leqslant x_1 + x_1' < 2x_2$, we deduce~$x_1' = - x_1 $ and hence~$x_1 = x_1' =0$
\proofend

\begin{proposition}
\label{prop:s2s2mon}
Let~$0 \leqslant x_1 \leqslant x_2$. Then the Hamiltonian monodromy group of the toric fibre~$T(x_1,x_2)\subset S^2 \times S^2$ in the monotone~$S^2 \times S^2$ is given by
\begin{equation}
\ch_{T(x_1,x_2)} = 
\begin{cases}
	\left \langle
	\begin{pmatrix}
		-1 & 0 \\
		0 & 1  
	\end{pmatrix}
 	\right \rangle \cong \Z_2, 
 	& x_1 = 0, x_2 \neq 0; \\[3ex]
	\left \langle
	\begin{pmatrix}
		0 & 1 \\
		1 & 0  
	\end{pmatrix}
 	\right \rangle \cong \Z_2,
 	& x_1 = x_2 \neq 0; \\[3ex]
	\left \langle
	\begin{pmatrix}
		1 & 0 \\
		0 & -1  
	\end{pmatrix},
	\begin{pmatrix}
		-1 & 0 \\
		0 & 1  
	\end{pmatrix}
 	\right \rangle \cong \Z_2 \times \Z_2,
 	& x_1 = x_2 = 0; \\[3ex]
\end{cases}
\end{equation}
and by~$\ch_{T(x_1,x_2)} = \{1\}$ in all other cases. 
\end{proposition}

Note that any other toric fibre is Hamiltonian isotopic to a fibre with~$0 \leqslant x_1 \leqslant x_2$, meaning that its Hamiltonian monodromy group is conjugate to one of the above. Thus the only isomorphism types of groups which appear are~$\Z_2, \Z_2 \times \Z_2$ and the trivial group. The astute reader may have wondered why~$\ch_{(0,0)}$ is not the full symmetry group of~$\Delta = [-1,1] \times [-1,1]$. This comes from the fact that some of these symmetries act non-trivially on~$H_2(S^2 \times S^2)$ (by exchanging the obvious generators) and thus they can be realized by symplectomorphisms, but not by Hamiltonian diffeomorphisms. We refer to~\cite{AugSmiWor22}, where the monodromy group generated by symplectomorphisms is determined for monotone toric fibres.
\smallskip

\proofof{Proposition~\ref{prop:s2s2mon}}
Again, the construction side can be obtained by symmetric probes and Theorem~\ref{thm:symmprobe}. For the obstruction side, let~$T(x_1,x_2)$ be a toric fibre of~$S^2\times S^2$ and~$\phi \in \Ham(S^2 \times S^2)$ a Hamiltonian diffeomorphism mapping this fibre to itself. We again analyse the map~$\phi_* \in \Aut H_2(X,T(x_1,x_2))$ and use the fact that~$\phi_*$ determines~$(\phi\vert_{T(x_1,x_2)})_*$.
Let us start with the case of the monotone fibre~$T(0,0) \subset S^2 \times S^2$, which is also a special case of~\cite[Theorem 2]{AugSmiWor22}. In this case, the distiguished classes are~$\cd(0,0) = \{D_1,D_2,D_3,D_4\}$. Therefore Theorem~\ref{thm:ambientmonodromy} implies that the ambient monodromy is a permutation of these classes. Since~$D_1 + D_3, D_2 + D_4 \in H_2(X)$, these two classes must be preserved under~$\phi_*$ which implies the claim. In the case~$x_1 = x_2 \neq 0$, the distinguished classes are~$\cd(x_1,x_1) = \{D_1,D_2\}$ and hence only permutations of~$D_1$ and~$D_2$ are permitted by Theorem~\ref{thm:ambientmonodromy}. Now let~$0 \leqslant x_1 < x_2$. Then the set of distinguished classes is~$\cd(x_1,x_2) = \{D_2\}$ and hence the ambient monodromy map takes~$D_2$ to~$D_2$. We set~$x_1 = x_1'$ in~\eqref{eq:xdiff1} and~\eqref{eq:xdiff2}. In the first case, we find that~$a_1 + a_2 = 1$ and a computation using the expressions for~$a_3$ and~$a_4$ from the proof of Proposition~\ref{prop:s2s2class}, this yields that the monodromy is trivial. In the second case, we find that~$x_1 = 0$ and~$a_1 + a_2 = 0$ and a similar computation shows that the monodromy maps~$e_1 \mapsto -e_1$ in that case. We conclude that the monodromy group is trivial whenever~$x_1 \neq 0$ and that it is generated by the map~$e_1 \mapsto -e_1, e_2 \mapsto e_2$ in case~$x_1 = 0$.
\proofend

\begin{remark}
We point out that in case~$S^2 \times S^2$ is equipped with a non-monotone symplectic form, the classification as well as the Hamiltonian monodromy is drastically different. Indeed, some equivalence classes~$\fH_x$ of fibres have accumulation points in~$\Delta$ and some fibres have infinite Hamiltonian monodromy groups. We refer to~\cite{BreKim23} for details.
\end{remark}

\subsection{The case of~$X = \CP^2$}
\label{ssec:cp2}

Let~$\CP^2$ be equipped with the symplectic form and moment polytope as in Example~\ref{ex:cp2}. We give the classification of toric fibres and the Hamiltonian monodromy groups without proof since the proofs are the same as for~$S^2 \times S^2$. The classification of toric fibres was first given in~\cite[Proposition 7.1]{SheTonVia19}. Note that all equivalences and Hamiltonian monodromies in the case of the monotone~$S^2 \times S^2$ are induced by symmetries of the moment polytope. The same holds in the case of~$\CP^2$.

\begin{proposition}
Toric fibres~$T(x),T(y) \subset \CP^2$ are equivalent if and only if~$x$ can be mapped to~$y$ by an integral symmetry of~$\Delta$. Similarly, the Hamiltonian monodromy group~$\ch_{T(x)}$ consists of transformations induced by integral symmetries of~$\Delta$ fixing the point~$x$.
\end{proposition}

\subsection{The case of~$X = \C \times S^2$}
\label{ssec:excs2}

Let~$\C \times S^2$ be equipped with the symplectic form~$\omega = \omega_{\C} \oplus \omega_{S^2}$. We normalize the moment map such that its moment polytope is given by~$\Delta = \R_{\geqslant -1} \times [-1,1]$. There are symmetric probes with directional vector~$e_2^* + k e_1^*$ for every~$k \in \Z$. The probe with~$k=0$ is admissible everywhere. For~$k = \pm 1$, the probes are admissible everywhere except when they hit a vertex of~$\Delta$. The symmetric probes with~$k \notin \{-1,0,1\}$ are admissible whenever they hit both facets~$\R_{\geqslant -1} \times \{1\}$ and~$\R_{\geqslant -1} \times \{-1\}$. As we shall see the latter types of symmetric probes are obsolete as all results can be proven using only those with~$k \in \{-1,0,1\}$.

\begin{proposition}
\label{prop:classcs2}
The classification of fibres in~$(X,\omega)$ is as follows. For~$x= (x_1,0) \in \Int\Delta$ with~$x_1 \geqslant 0$, we have~$\fH_x = \{x\}$, i.e.\ the corresponding toric fibre is not equivalent to any other fibres. For~$x = (x_1,\pm x_1)$ with~$x_1 < 0$, we have~$\fH_x = \{(x_1,x_1),(x_1,-x_1)\}$. For~$x = (0,x_2)$ with~$x_2 > 0$, we have
\begin{equation}
	\fH_x = \{ (2nx_2, \pm x_2) \, \vert \, n \in \N \} \cup \{(-x_2,0)\}.
\end{equation}
For~$x = (x_1,\pm x_1)$ with~$x_1 > 0$, we have
\begin{equation}
	\fH_x = \{ ((2n+1)x_1, \pm x_1) \, \vert \, n \in \N \}.
\end{equation}
For~$x = (x_1,x_2)$ with~$0 < x_1 < x_2$, we have
\begin{equation}
	\fH_x = \{ (\pm x_1 + 2nx_2, \pm x_2) \, \vert \, n \in \N \} \cup \{(-x_2,\pm x_1)\}.
\end{equation}
All~$y \in \Int\Delta$ are in one of the above~$\fH_x$. 
\end{proposition}

\begin{figure}
	\begin{center}
	\begin{tikzpicture}
		\node[inner sep=0pt] at (0,0)
    			{\includegraphics[trim={3cm 6cm 2.5cm 1cm},clip,scale=0.8]						{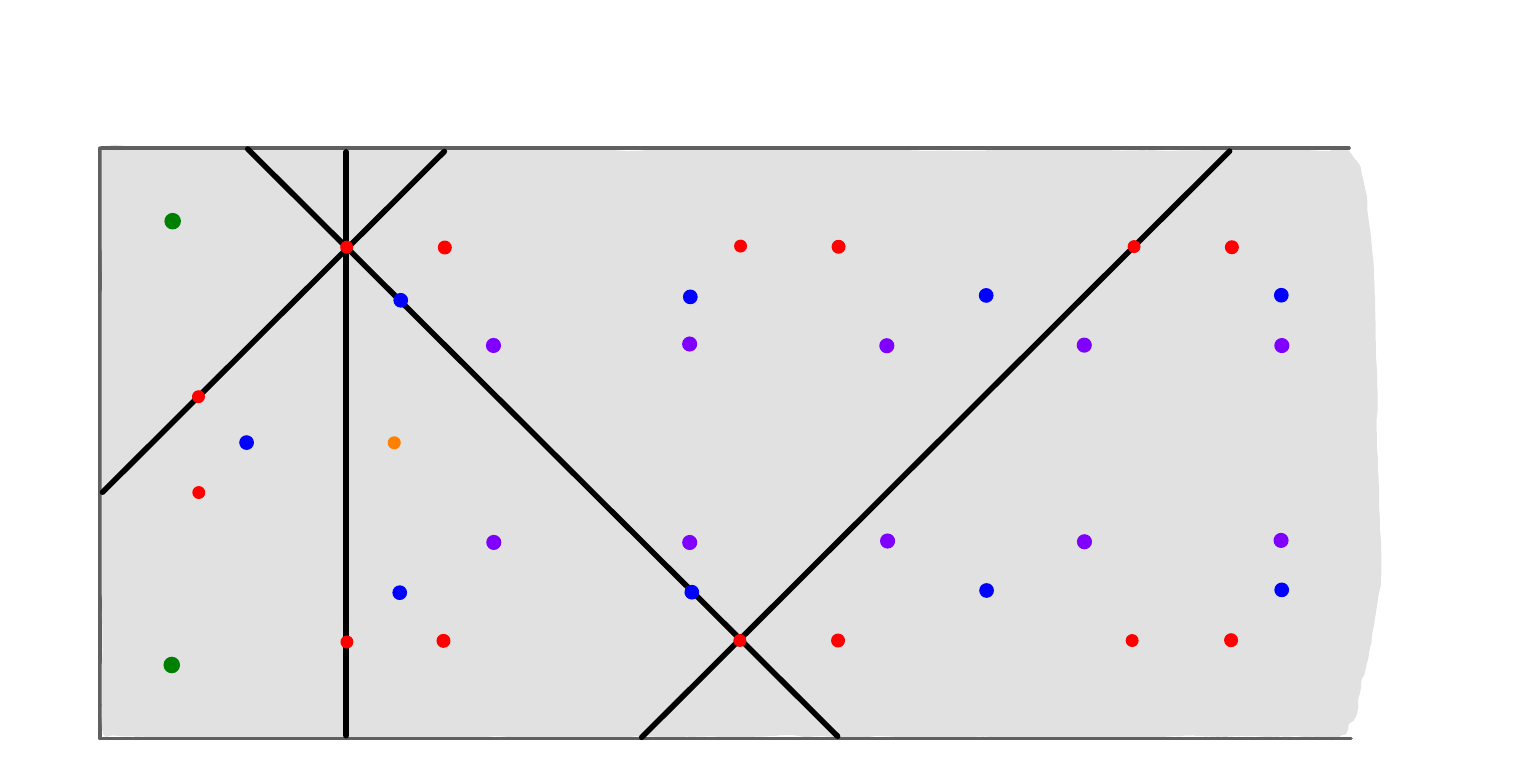}};
	\end{tikzpicture}
	\caption{Some symmetric probes in~$\C \times S^2$. Points of the same colour denote equivalent fibres.}
	\label{fig:7}
	\end{center}
\end{figure}

\proof
For the construction of the equivalences, we use concatenations of the probes with directional vectors~$e_2^* + k e_1^*$ for~$k \in \{-1,0,1\}$, as we discuss below on a case by case basis. For the obstruction side, note that we cannot directly apply Theorems~\ref{thm:mainB} and~\ref{thm:ambientmonodromy}, since~$X$ is non-compact. However,~$X$ is of reduction type, see Definition~\ref{def:redtype}. Indeed, the toric reduction of~$\C^2$ to~$S^2$ coming from the Delzant construction yields a toric reduction of~$\C \times \C^2$ to~$\C \times S^2$. Therefore, the results of Theorems~\ref{thm:mainB} and~\ref{thm:ambientmonodromy} still apply. 

By the first Chekanov invariant from Theorem~\ref{thm:mainB}, the polytope~$\Delta$ decomposes into subsets~$\fD_x$ of constant distance~$0 < d(x) \leqslant 1$ to the boundary~$\pp\Delta$. First, let~$0 < d(x) < 1$. The toric fibres of the type~$(x_1,\pm x_1)$ with~$x_1 < 0$ are the only ones having~$\#_d = 2$, which distinguishes them from all others. Note that for any other~$x \in \Delta$ with~$d(x) < 1$, the torus~$T(x)$ is equivalent by symmetric probes to exactly one torus on the segment~$[0,x_2] \times \{x_2\}$ with~$x_2 = 1-d(x)$. It is easy to see that by probes with directional vectors~$e_2 + e_1$ and~$e_2 - e_1$, any fibre is equivalent to one on the segment~$[-x_2,x_2] \times \{x_2\}$. Now note that fibres on this segment come in equivalent pairs as can be seen by noting that~$(x_1,x_2)$ is equivalent to~$(-x_2,-x_1)$ which is equivalent (by a vertical probe) to~$(-x_2,x_1)$ which in turn is equivalent to~$(-x_1,x_2)$. Thus the problem of classifying fibres with~$d(x)<1$ boils down to classifying fibres on the segment~$[0,x_2] \times \{x_2\}$. 

\textbf{Claim:} If $T(x_1,x_2) \cong T(x_1',x_2)$ for~$x_1,x_1' \in [0,x_2]$ and~$x_2 = 1 - d(x)$, then~$x_1 = x_1'$.

\ni
To prove the claim, we follow the same strategy as in the proof of Proposition~\ref{prop:s2s2class}. Suppose there is~$\phi \in \Ham(X,\omega)$ mapping~$T(x_1,x_2)$ to~$T(x_1',x_2)$ and let~$\phi_*$ be the ambient monodromy induced by this map. Let~$D_1,D_2,D_3$ be the canonical basis of~$H_2(X,T(x_1,x_2))$ where~$D_1$ is the disk corresponding to the facet~$\{-1\} \times [-1,1]$ and the remaining ones ordered in the anti-clockwise direction. Let~$D_1',D_2',D_3' \in H_2(X,T(x_1',x_2))$ be the disks obtained by the same convention. The distinguished classes are~$\cd(x_1,x_2) = \{D_3\}$ and~$\cd(x_1',x_2) = \{D_3'\}$ meaning that~$\phi_*D_3 = D_3'$. Set
\begin{equation}
	\label{eq:phiD}
	\phi_* D_1 = a_1D_1' + a_2D_2' + a_3D_3', \quad
	a_i \in \Z.
\end{equation}
By the invariance of the Maslov class, we obtain~$a_1 + a_2 + a_3 = 1$. Since the induced map~$(\phi\vert_{T(x_1,x_2)})_*$ is invertible, we deduce that~$\det(\pp\phi_*D_1,\pp\phi_*D_3) = \pm 1$ which yields~$a_1 = \pm 1$. Preservation of symplectic area,~$\int_{D_1} \omega = \int_{\phi_*D_1} \omega$ yields
\begin{equation}
	1 + x_1 = a_1(1 + x_1') + a_2(1 + x_2) + a_3(1 - x_2).
\end{equation}
In the case~$a_1 = -1$, we find
\begin{equation}
	\label{eq:xdiff3}
	x_1 + x_1' = 2x_2(a_2 - 1), 
\end{equation}
from which we deduce that~$x_1 = x_1'$. Similarly, for~$a_1 = 1$, we find
\begin{equation}
	\label{eq:xdiff4}
	x_1 - x_1' = 2x_2a_2, 
\end{equation}
which implies the same conclusion and thus proves the claim.

Let us now turn to the case~$d(x)=1$, i.e.\ tori of the form~$T(x_1,0)$ with~$x_1 \geqslant 0$. Note that~$T(0,0)$ is the only monotone fibre and thus not equivalent to any other fibre. We will show that the same remains true for~$x_1 > 0$. Indeed, if~$T(x_1,0)$ and~$T(x_1',0)$ were equivalent, then the same arguments as above apply to the ambient monodromy~$\phi_*$ except that now~$\cd(x_1,0) = \{D_2,D_3\}$. Equations~\eqref{eq:xdiff3} and~\eqref{eq:xdiff4} for~$x_2 = 0$ imply the claim. Equivalently, one can use~\cite[Theorem 1.1]{FOOO13} to find that~$e(X,T(x_1,0)) = 1 + x_1$, meaning that these fibres are also distinguished by their displacement energy. 
\proofend

\begin{proposition}
Let~$(x_1,x_2) \in \Int \Delta$. The Hamiltonian monodromy group of the toric fibre~$T(x_1,x_2) \subset \C \times S^2$ is given by
\begin{equation}
\label{eq:moncs2}
\ch_{T(x_1,x_2)} = 
\begin{cases}
	\left \langle
	\begin{pmatrix}
		0 & -1 \\
		-1 & 0  
	\end{pmatrix}
 	\right \rangle \cong \Z_2,
 	& x_1 = - x_2, x_2 > 0; \\[3ex]
	\left \langle
	\begin{pmatrix}
		1 & 0 \\
		0 & -1  
	\end{pmatrix}
 	\right \rangle \cong \Z_2,
 	& x_1 \leqslant 0 , x_2 = 0; \\[3ex]
 	\left \langle
	\begin{pmatrix}
		1 & 0 \\
		0 & -1  
	\end{pmatrix},
	\begin{pmatrix}
		1 & 0 \\
		2 & -1  
	\end{pmatrix}
 	\right \rangle \cong \Z_2 \ltimes \Z,
 	& x_1 > 0, x_2 = 0;
\end{cases}
\end{equation}
and is given by a group conjugated to the above in case~$T(x_1,x_2)$ is equivalent to one of the above cases according to Proposition~\ref{prop:classcs2}, and by~$\ch_{T(x_1,x_2)} = \{1\}$ in all other cases. 
\end{proposition}

\proof
In case~$x_1 = - x_2, x_2 > 0$, we have~$\cd(x_1,x_2) = \{D_1,D_3\}$ and thus any monodromy matrix has to permute~$\pp D_1 = e_1$ and~$\pp D_3 = -e_2$. This permutation is realized by a symmetric probe with directional vector~$e_1^* + e_2^*$. For the case~$x_1 \leqslant 0 , x_2 = 0$, note that the vertical symmetric probe realizes the claimed monodromy. In terms of obstructions, note that~$\cd(0,0) = \{D_1,D_2,D_3\}$, which yields the claim for~$T(0,0)$. In case~$x_1 < 0$, note that~$T(x_1,0)$ is Hamiltonian isotopic to~$T(0,x_1)$ to which we can apply the equations~\eqref{eq:xdiff3} and~\eqref{eq:xdiff4} to find that the only possible monodromy for~$T(0,x_1)$ is~$e_1 \mapsto -e_1, e_2 \mapsto e_2$. Under the conjugation induced by the equivalence of~$T(x_1,0)$ and~$T(0,x_1)$, this yields the answer. Now let~$x_1 \in (0,x_2]$ and~$x_2 > 0$. Then the Hamiltonian monodromy group is trivial by equations~\eqref{eq:xdiff3} and~\eqref{eq:xdiff4}. Now let us turn to the case of the infinite monodromy groups for the fibres~$T(x_1,0)$ with~$x_1 > 0$. The two generators given in~\eqref{eq:moncs2} correspond to the vertical probe and the probe with directional vector~$-e_1^* + e_2^*$. Since~$\cd(x_1,0) = \{D_2,D_3\}$, we distiguish the cases~$D_2 \mapsto D_2, D_3 \mapsto D_3$ and~$D_2 \mapsto D_3, D_3 \mapsto D_2$. Let us first restrict our attention to the former case. We use~\eqref{eq:phiD} with~$D_i = D_i'$. Recall that~$a_1 = \pm 1$. Since~\eqref{eq:xdiff3} cannot be satisfied for~$x_1 = x_1', x_2 = 0$, we deduce that~$a_1 = 1$. A computation shows that
\begin{equation}
	\label{eq:monocs2}
	(\phi\vert_{T(x_1,0)})_* e_1 
	= \pp \phi_* D_1
	= e_1 + 2a_2 e_2,
\end{equation}
which proves the claim in case~$\det(\phi\vert_{T(x_1,0)})_* = 1$ (this corresponds to the powers of the product of the two generators given in~\eqref{eq:moncs2}). The case~$D_2 \mapsto D_3, D_3 \mapsto D_2$ is completely analoguous.
\proofend

\subsection{The case of~$X = \C^2 \times T^*S^1$}
\label{ssec:c2ts1}

Let~$X = \C^2 \times T^*S^1$ be equipped with the exact symplectic form~$\omega = d\lambda = \omega_{\C^2} \oplus \omega_{T^*S^1} = d\lambda_{\C^2} \oplus d\lambda_{T^*S^1}$ and the product toric structure with moment polytope~$\Delta = \R_{\geqslant 0}^2 \times \R$. Note that~$X$ is not a toric reduction of~$\C^N$ and hence we cannot apply the techniques from~\S\ref{ssec:chekanovinv} and~\S\ref{ssec:hammon} as in the previous examples. Instead, we rely on~\S\ref{ssec:displacementen}, meaning that we need to compute the displacement energy of toric fibres. 

\begin{lemma}
\label{lem:c2ts1}
The displacement energy of toric fibres is given by
\begin{equation}
	e(X,T(x_1,x_2,x_3))
	=
	\min\{x_1,x_2\}.
\end{equation}
In particular, equality~\eqref{eq:detoric} (and thus also Assumption~\ref{ass:detoric}) holds for all toric fibres in~$X$.
\end{lemma}

\proof
The upper bound is obvious, either by using probes or by using the fact that $e(\C^2,T(x_1,x_2)) = \min\{x_1,x_2\}$. For the lower bound, we use Chekanov's theorem~\cite{Che98}, which we briefly recall here. Let~$L \subset X$ be a compact Lagrangian submanifold of a tame symplectic manifold and let~$J \in \cj(X,\omega)$ be a tame almost complex structure. Furthermore, denote by~$\sigma(X,L;J)$ the infimum of symplectic areas over all non-constant~$J$-holomorphic disks with boundary on~$L$. If this set is empty, set~$\sigma(X,L;J) = \infty$. If the set is not empty, we obtain a strictly positive value which is attained by Gromov compactness. The quantity~$\sigma(X;J)$ is defined similarly for~$J$-holomorphic spheres in~$X$. Then Chekanov's theorem gives the lower bound,
\begin{equation}
	e(X,L) \geqslant \min \{ \sigma(X;J) , \sigma(X,L;J) \}.
\end{equation}
Now let~$X = \C^2 \times T^*S^1$ and~$L = T(x_1,x_2,x_3)$ a toric fibre. Note that~$X$ is aspherical, and thus~$\sigma(X;J)= \infty$. Let~$J_0 \in \cj(X,\omega)$ be the complex structure obtained from the identification~$X = \C^2 \times \C^{\times}$. There are two obvious families of~$J_0$-holomorphic disks, 
\begin{eqnarray}
	u_{\alpha_1,\alpha_2}(z) &=& \left(z,\sqrt{\frac{x_2}{\pi}}e^{i\alpha_1} , e^{x_3 + i\alpha_2} \right), \\
	v_{\alpha_1,\alpha_2}(z) &=& \left(\sqrt{\frac{x_1}{\pi}}e^{i\alpha_1}, z , e^{x_3 + i\alpha_2} \right),\quad \alpha_1, \alpha_2 \in S^1.
\end{eqnarray}
These disks have area~$\int u_{\alpha_1,\alpha_2}^*\omega = x_2$ and~$\int v_{\alpha_1,\alpha_2}^*\omega = x_1$ for all~$\alpha \in S^1$. We show that the minimal one among these two disks realizes the minimum~$\sigma(X,L;J_0)$. For a similar argument, see~\cite[Lemma 4]{Aur15}. Now let~$u \colon (D,\pp D) \rightarrow (X,L)$ be a non-trivial~$J_0$-holomorphic disk. First note that the map~$p_2 \circ u$, where~$p_2 \colon X \rightarrow T^*S^1 \cong \C^*$ is the projection, is constant. This follows from the maximum principle. Indeed, by the maximum principle this map takes values in the unit disk. Since~$0$ is not contained in its image, we can precompose it with~$z \mapsto \frac{1}{z}$, to see that it actually takes values in the unit circle. Since it is holomorphic, it is actually constant. By considering~$p_1 \circ u$, where~$p_1 \colon (X,J_0) \rightarrow (\C^2,i \oplus i)$ is the natural projection, it is sufficient to understand holomorphic disks in~$\C^2$ with boundary on the product torus~$T(x_1,x_2)$. The group~$\pi_2(\C^2,T(x_1,x_2))$ is generated by the two coordinate disks~$D_1,D_2$. We have~$[p_1 \circ u] = k_1D_1 + k_2D_2$, where~$k_1,k_2$ are the algebraic intersection numbers with coordinate axes. By positivity of intersections, we deduce~$k_1,k_2 \geqslant 0$. Since~$\int_{D_1} \omega_{\C^2} = x_1$ and~$\int_{D_2} \omega_{\C^2} = x_2$, we obtain
\begin{equation}
	\sigma(X,L;J_0)
	\geqslant \min\{ k_1x_1 + k_2x_2 \,\vert\, k_1,k_2 \geqslant 0, k_1k_2 \neq 0 \}
	= \min\{x_1,x_2\}.
\end{equation}
This minimum is realized by the disks~$u_{\alpha_1,\alpha_2}$ or $v_{\alpha_1,\alpha_2}$.
\proofend

Using Theorem~\ref{thm:alternative}, we can now classify toric fibres and determine their Hamiltonian monodromy groups.

\begin{proposition}
\label{prop:classc2ts1}
The classification of toric fibres~$T(x) = T(x_1,x_2,x_3)$ in~$X = \C^2 \times T^*S^1$ is given by
\begin{equation}
	\fH_x = \{(x_1,x_2,x_3 + k(x_2 - x_1)) , (x_2,x_1,x_3 + k(x_2 - x_1)) \,\vert\, k \in \Z\}.
\end{equation}
Furthermore, all Hamiltonian monodromy groups are trivial except when~$x_1 = x_2$, in which case
\begin{equation}
	\label{eq:c2ts1mon}
	\ch_{T(x_1,x_1,x_3)}
	=
	\left\langle 
		\begin{pmatrix}
		1 & 0 & 1 \\
		0 & 1 & -1 \\
		0 & 0 & 1
		\end{pmatrix} ,
		\begin{pmatrix}
		0 & 1 & 0 \\
		1 & 0 & 0 \\
		0 & 0 & 1
		\end{pmatrix}
	\right\rangle.
\end{equation}
\end{proposition}

\proof 
The equivalences are easy to construct using probes with direction~$e_1^* - e_2^* + k e_3^*$ with~$k \in \Z$. Let~$\phi$ be a Hamiltonian diffeomorphism mapping~$T(x)$ to~$T(x')$. First note that the long exact sequence for relative homology looks quite different than in the compact toric case. Indeed, we obtain
\begin{equation}
	0 
	\rightarrow H_2(X,T(x)) 
	\rightarrow H_1(T(x))
	\rightarrow H_1(X)
	\rightarrow 0,
\end{equation}
and~$H_1(T(x)) \cong H_2(X,T(x)) \oplus H_1(X) \cong \Z^2 \oplus \Z$.  Then the induced map~$(\phi\vert_{T(x)})_*$ on the first homology is of the form 
\begin{equation}
	(\phi\vert_{T(x)})_*
	=
	\begin{pmatrix}
		a_1 & a_2 & b_1 \\
		a_3 & a_4 & b_2 \\
		0 & 0 & 1
	\end{pmatrix} , \quad a_i,b_j \in \Z.
\end{equation}
This follows from the fact that~$\phi$ induces the identity on~$H_1(X)$. First suppose that~$x_1 = x_2$. Then, by Theorem~\ref{thm:alternative} and Lemma~\ref{lem:c2ts1}, we obtain~$x_1' = x_2' = x_1$ and~$(\phi\vert_{T(x)})_*$ either swaps the first two coordinates or acts by the identity on them. By preservation of the Maslov index, we obtain~$b_2 = -b_1$. Note that all monodromies of these tori can be realized by symmetric probes of direction~$e_1^* - e_2^* + k e_3^*$, proving~\eqref{eq:c2ts1mon}. To prove that~$x_3 = x_3'$, compute
\begin{equation}
	\label{eq:c2ts1}
	x_3 	
	= \int_{e_3} \lambda \vert_{T(x)}
	= \int_{(\phi\vert_{T(x)})_*e_3} \lambda \vert_{T(x')}
	= b_1 (x_1' - x_2') + x_3,
\end{equation}
proving the claim. In the case~$x_1 \neq x_2$, suppose without loss of generality that~$x_1 < x_2$ and~$x_1' < x_2'$. By Theorem~\ref{thm:alternative} and Lemma~\ref{lem:c2ts1}, we obtain~$x_1 = x_1'$ and~$a_1 = 1, a_3 = 0$ and~$a_4=1, a_2 = 0$ or~$a_4 = -1, a_2 = -1$. The latter case is actually impossible. Indeed, if~$a_4 = -1, a_2 = -1$, then
\begin{equation}
	x_2 
	= \int_{e_2} \lambda \vert_{T(x)}
	= \int_{(\phi\vert_{T(x)})_*e_2} \lambda \vert_{T(x')}
	= 2 x'_1 - x'_2, 
\end{equation}
contradicting~$x_1 < x_2, x_1' < x_2'$. The rest of the proof is as in the case~$x_1 = x_2$ and the claim~$b_1 = 0$ follows from~\eqref{eq:c2ts1}.
\proofend 

\subsection{The case of~$X = T^*S^1 \times S^2$}
\label{ssec:ts1s2}

Let~$X = T^*S^1 \times S^2$ be equipped with the product symplectic form~$\omega = \omega_{T^*S^1} \oplus \omega_{S^2}$. The moment polytope is~$\Delta = \R \times [-1,1]$. Note that~$X$ is not a toric reduction of any~$\C^N$, but it is a toric reduction of~$\C^2 \times T^*S^1$, meaning that we can use Proposition~\ref{prop:classc2ts1} together with the lifting trick. 

\begin{proposition}
The classification of toric fibres~$T(x) = T(x_1,x_2)$ in~$X = T^*S^1 \times S^2$ is given by
\begin{equation}
	\fH_x 
	= 
	\{ (x_1 + 2k x_2,\pm x_2) \, \vert \, k \in \Z \}.
\end{equation}
Furthermore, all Hamiltonian monodromy groups are trivial, except for~$x_2 = 0$, in which case, 
\begin{equation}
	\ch_{T(x_1,0)}
	=
	\left\{ \left.
		\begin{pmatrix}
		1 & 0 \\
		2k & \pm 1
		\end{pmatrix}
		\right\vert 
		k \in \Z
	\right\}.
\end{equation}
\end{proposition}

\proof
Again, all constructions immediately follow from symmetric probes. For the obstructions, we view~$X$ as a toric reduction of~$\C^2 \times T^*S^1$. Perform toric reduction on~$\C^2 \times T^*S^1$ with respect to the plane~$V = \{ x_1 + x_2 = 1 \} \subset \R^3$ to obtain~$X$. This corresponds to the Hamiltonian~$H = \pi \vert z_1 \vert^2 + \pi \vert z_2 \vert^2$. The classification of toric fibres follows immediately from Propostion~\ref{prop:classc2ts1}, Proposition~\ref{prop:toricfibresred} and the lifting trick, as in the proof of Theorem~\ref{thm:mainB}. The obstructions to monodromy similarly follow from Proposition~\ref{prop:classc2ts1} and the lifting trick, as in the proof of Theorem~\ref{thm:ambientmonodromy}. \proofend

\subsection{Chekanov's classification revisited}
\label{ssec:chekrev}

In this subsection, we prove Conjecture~\ref{conj:main} in the case of~$\C^n$. The classification of product tori goes back to Chekanov, see Theorem~\ref{thm:chekanov}, meaning that we only need to prove that all equivalences of toric fibres can be realized by iterated symmetric probes. 

\begin{theorem}
\label{thm:cnsymmprobes}
Product tori~$T(x),T(y) \subset \C^n$ are Hamiltonian isotopic if and only if they are equivalent by a sequence of symmetric probes.
\end{theorem}

Before proving this result, let us revisit Chekanov's classification. In~$\C^2$, it states that,
\begin{equation}
	\fH_x = \{(x_1,x_2),(x_2,x_1)\}, \quad (x_1,x_2) \in \R^2_{> 0}.
\end{equation}
In~$\C^n$ with~$n \geqslant 3$, however, the situation is much richer. Note for example that all tori~$T(1,2,k)$ with~$k \in \N_{\geqslant 2}$ are Hamiltonian isotopic, since their Chekanov invariants agree. The set~$\fH_x$ even has accumulation points in many cases, see Corollary~\ref{cor:accumulation}. To discuss this further, we slightly reformulate Chekanov's invariants. Since coordinate permutations can be realized by Hamiltonian isotopies, we may assume that~$T(x)$ is given under the form
\begin{equation}
\label{eq:productnormal}
	T(x) = T(\underbrace{\underline{x},\dots,\underline{x}}_{\#_d(x)},\underline{x}+\xh_1 ,\dots,\underline{x}+\xh_s),
\end{equation}
for~$\xh_i > 0$ and~$s = n - d(x)$. We call~$\xh = (\xh_1, \ldots, \xh_s) \in \R_{>0}^s$ the \emph{reduced vector} associated to~$x$. 

Theorem~\ref{thm:chekanov} can be reformulated as follows. 

\begin{corollary}
\label{cor:chekref}
Product tori~$T(x),T(y) \subset \C^n$ are Hamiltonian isotopic if and only if~$d(x) = d(y)$,~$\#_d(x) = \#_d(y)$ and there is~$M \in \GL(s;\Z)$ such that~$M\xh = \yh$.
\end{corollary}

\proof
Note that the~$\xh_i$ are exactly the non-trivial generators of the lattice~$\Gamma(x)$,
\begin{equation}
	\Gamma(x) 
	= \Z\langle \xh \rangle
	= \{ k_1 \xh_1 + \ldots + k_s \xh_s \, \vert \, k_i \in \Z \} \subset \R.
\end{equation}
Furthemore,~$\Z\langle \xh \rangle$ is a complete invariant for~$\GL(s;\Z)$-orbits. See for example Cabrer--Mundici~\cite[Proposition 1]{CabMun16}. 
\proofend

This allows us to gain a good qualitative understanding of~$\fH_x$. 

\begin{corollary}
The inclusion 
\begin{equation}
\label{eq:fhinclusion}
	\fH_x \subset \{y \in \R^n_{>0} \, \vert \, d(y) = d(x), \, \#_d(y) = \#_d(x)\}
\end{equation}
is dense if and only if~$\rank \Gamma(x) \geqslant 2$.
\end{corollary} 

\proof
Let~$\xh \in \R_{>0}^s$ be the reduced vector as in~\eqref{eq:productnormal}. It follows from Corollary~\ref{cor:chekref} that the inclusion~\eqref{eq:fhinclusion} is dense if and only if the~$\GL(s;\Z)$-orbit of~$\xh$ is dense in~$\R_{>0}^s$. The latter is equivalent to~$\rank \Gamma(x) \geqslant 2$ by a classical theorem of Dani~\cite[Theorem 17]{Dan75}, see also~\cite{CabMun16}.
\proofend

Let us now have a look at the discrete case, i.e.\ the case where~$\rank(\Gamma(x))=1$. In that case, the reduced vector~$\xh \in \R^s_{>0}$ is a real multiple of a lattice vector,~$\xh = \intl(\xh) k$ with~$k \in \Z^s$ a primitive vector and where~$\intl(\xh) > 0$ denotes the integral affine length.

\begin{corollary}
\label{cor:rankone}
The product tori~$T(x),T(y)\subset \C^n$ with~$d(x)=d(y)$,~$\#_d(x)=\#_d(y)$ and~$\rank(\Gamma(x))=\rank(\Gamma(y)) = 1$ are Hamiltonian isotopic if and only if their reduced vectors have the same integral affine length,
$\intl(\xh)=\intl(\yh)$.
\end{corollary}

\proof
Write~$\xh = \intl(\xh) k$,~$\yh = \intl(\yh) k'$ and note that~$\GL(s;\Z)$ acts transitively on the set of primitive lattice vectors and preserves integral affine length. Thus the claim follows from Corollary~\ref{cor:chekref}.
\proofend

Let us now turn to the proof of Theorem~\ref{thm:cnsymmprobes}. The following lemma is key.

\begin{lemma}
\label{lem:elswap}
Let~$x=(x_1,x_2,x_3) \in \R^3_{> 0}$ with~$x_1 < x_2,x_3$. Then there is a symmetric probe showing
\begin{equation}
	\label{eq:elswap}
	T(x_1,x_2,x_3) \cong T(x_3,x_2 + x_3 - x_1, x_1) \subset \C^3.
\end{equation}
\end{lemma}

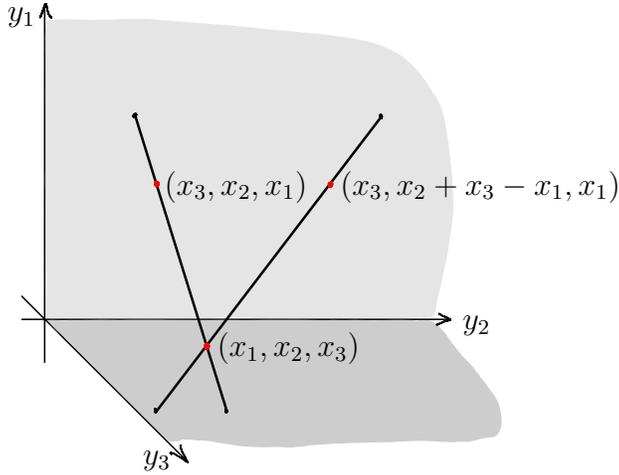
\begin{figure}
	\begin{center}
	\begin{tikzpicture}
		\node[inner sep=0pt] at (0,0)
    			{\includegraphics[trim={4cm 3cm 6cm 1.5cm},clip,scale=0.6]						{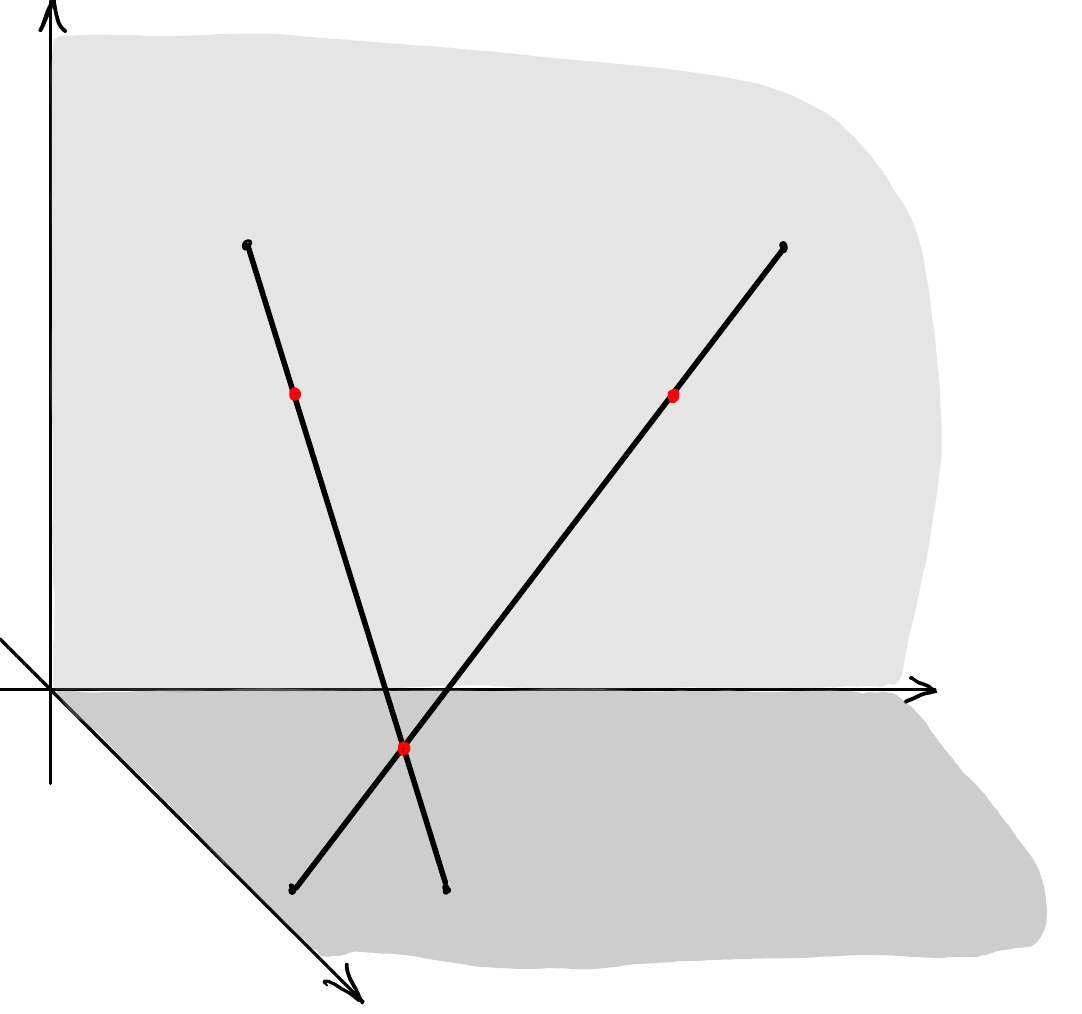}};
    		\node at (0.25,-1.5){$(x_1,x_2,x_3)$};
    		\node at (-0.45,0.65){$(x_3,x_2,x_1)$};
    		\node at (2.8,0.65){$(x_3,x_2+x_3-x_1,x_1)$};
    		\node at (-3.3,2.9){$y_1$};
    		\node at (2.75,-1.2){$y_2$};
    		\node at (-1.5,-3){$y_3$};
	\end{tikzpicture}
	\caption{Two symmetric probes in~$\C^3$, one realizes a coordinate permutation, the other is used in the proof of Lemma~\ref{lem:elswap}}
	\label{fig:5}
	\end{center}
\end{figure}

\proof 
The directional vector~$\eta = e_1^* +e_2^* - e_3^*$ defines a probe~$\sigma = \R^3_{\geqslant 0} \cap \{ x + t\eta \, \vert \, t \in \R \}$ which realizes the equivalence~\eqref{eq:elswap}, see Figure~\ref{fig:5}. Indeed, since~$x_1 < x_2$, the probe intersects the boundary~$\pp \R^3_{\geqslant 0}$ in the points~$(0,x_2-x_1,x_3+x_1)$ and~$(x_1+x_3,x_2+x_3,0)$ which lie in the interior of the facets~$\{y_1 = 0\}$ and~$\{y_3 = 0\}$ respectively. Since~$\langle \eta , e_1 \rangle = 1$ and~$\langle \eta , e_3 \rangle = -1$, both intersections are integrally transverse and hence the probe is admissible. Furthermore, since~$x_1 \leqslant x_2,x_3$, the points~$(x_1,x_2,x_3)$ and~$(x_3,x_2 + x_3 - x_1, x_1)$ both lie at integral distance~$x_1$ to the boundary and hence the corresponding fibres are Hamiltonian isotopic.
\proofend

\proofof{Theorem~\ref{thm:cnsymmprobes}}
Let~$T(x),T(y) \subset \C^n$ be product tori whose Chekanov invariants agree. First note that coordinate permutations 
\begin{equation}
	(x_1,\ldots,x_i,\ldots,x_j,\ldots,x_n) \mapsto
	(x_1,\ldots,x_j,\ldots,x_i,\ldots,x_n)
\end{equation}
can be realized by symmetric probes. Therefore, we may assume that both tori are given in the normal forms
\begin{equation}
	x=(\underline{x}, \ldots, \underline{x}, \underline{x} + \xh_1, \ldots, \underline{x} + \xh_s ), \quad
	y=(\underline{x}, \ldots, \underline{x}, \underline{x} + \yh_1, \ldots, \underline{x} + \yh_s ),
\end{equation}
where the reduced vectors~$\xh,\yh \in \R^s_{>0}$ are~$\GL(s;\Z)$-equivalent by Corollary~\ref{cor:chekref}. We thus need to prove that all~$\GL(s;\Z)$-transformations on the reduced vectors can be realized by symmetric probes. The group~$\GL(s;\Z)$ is generated by coordinate permutations together with the transformation
\begin{equation}	\label{eq:generator}
	(\xh_1,\xh_2,\ldots,\xh_s) \mapsto (\xh_1+\xh_2,\xh_2, \ldots, \xh_s),
\end{equation}
see for example~\cite{Tro62}. Coordinate permutations can be realized by symmetric probes as we have just seen. By Lemma~\ref{lem:elswap}, the generator~\eqref{eq:generator} can also be realized by a symmetric probe. Indeed, note that the reduced vectors associated to the tori~$T(x_1,x_2,x_3) \subset \C^3$ and~$T(x_3,x_2 + x_3 - x_1, x_1) \subset \C^3$ are given by~$(x_2-x_1,x_3-x_1) \in \R^2_{>0}$ and~$(x_2+x_3-2x_1, x_3 - x_1)\in \R^2_{>0}$ which corresponds to~$(\xh_1,\xh_2) \mapsto (\xh_1+\xh_2, \xh_2)$ in terms of the reduced vectors. Hence~\eqref{eq:generator} can be realized by a symmetric probe lying in an appropriately chosen coordinate subspace~$\C^3 \subset \C^n$. This proves the claim of the theorem.
\proofend

\begin{remark}
\label{rk:classquant}
Let us briefly discuss a quantitative version of Theorem~\ref{thm:cnsymmprobes}. More specifically, 
For every~$\varepsilon > 0$, we can find a Hamiltonian isotopy which has support in the ball~$B^6(x_1 + x_2 + 2x_3 + \varepsilon)$ realizing the equivalence in Lemma~\ref{lem:elswap}. Indeed, note that the closure of~$B^6(x_1 + x_2 + 2x_3)$ is the smallest closed ball containing the symmetric probe~$\sigma$ in the proof of Lemma~\ref{lem:elswap}, and the support of the Hamiltonian isotopy can be chosen to lie in an arbitrarily small neighbourhood of~$\sigma$. This yields a simple proof of~\cite[Lemma 4.1]{CheSch16}. Furthermore, by the same argument as in the proof of~\cite[Theorem 1.1 (ii)]{CheSch16}, this remark implies that for
\begin{equation}
	\label{eq:balladmissible}
	\sum_{i=1}^n x_i + d(x) < a, \quad
	\sum_{i=1}^n y_i + d(y) < a
\end{equation}
the product tori~$T(x),T(y)$ are Hamiltonian isotopic by iterated symmetric probes \emph{inside the ball}~$B^{2n}(a)$. In other words, given a ball~$B^{2n}(a)$, there is a region~$\cR(a) \subset B^{2n}(a)$ in its (open) moment polytope
\begin{equation}
	\cR(a) = \mu_0^{-1}\left\{ x \in \R^n_{\geqslant 0} \, \left\vert \, \sum_{i=1}^n x_i + d(x) < a \right.\right\}
\end{equation}
in which the classification of product tori in~$B^{2n}(a)$ coincides with the classification of product tori in all of~$\C^n$ and all symmetric probes producing the equivalences are contained in~$\mu_0(\cR(a))$. Together with Chekanov's classification theorem~\ref{thm:chekanov}, this shows Corollary~\ref{cor:accumulation}, which also follows from the methods in~\cite{CheSch16}. Let us point out that one cannot drop~\eqref{eq:balladmissible} as was shown in~\cite[Theorem 1.2]{CheSch16}. A reasonable guess would be that in the complement of~$\cR(a)$, only coordinate permutations are allowed, since these are the only symmetric probes admissible in that region. However, we do not know the classification of product tori in the ball outside of~$\cR(a)$. 
\end{remark}

\begin{corollary}
\label{cor:accumulation}
Let~$X$ be a toric manifold of dimension~$\geqslant 6$ whose moment polytope has at least one vertex. Then~$X$ contains toric fibres such that~$\fH_x$ has accumulation points. 
\end{corollary}

\subsection{In arbitrary toric manifolds}
\label{ssec:arbitrary}

The goal of this section is to illustrate that there are many symmetric probes in arbitrary toric manifolds. We focus on constructions near the boundary of the moment polytope~$\Delta$ using normal forms of Delzant polytopes. For example each vertex~$v \in \Delta$ yields an equivariant symplectic ball embedding~$B^{2n}(a) \rightarrow X$, for each~$a$ smaller than the integral length of the shortest edge adjacent to~$v$. Note that this ball embedding is unique up to coordinate permutations. Let us denote the corresponding subset by~$B_v(a) \subset X$. Furthermore, denote the image of a product torus~$T(x) \subset B^{2n}(a)$ under the equivariant embedding by~$T_v(x) \subset B_v(a)$. From Remark~\ref{rk:classquant}, we deduce that there is a region~$\cR_v(a) \subset B_v(a)$ in which the same probes are admissible as those that are admissible in~$\C^n$. Let us show the following result about the classification of toric fibres close to vertices. 

\begin{proposition}
\label{prop:vertexclass}
Let~$B^{2n}_v(a)$ and~$B^{2n}_{v'}(a)$ be balls at vertices~$v,v' \in \Delta$ with
\begin{equation}
	0 < a < \min\{ \intl(e) \, \vert \, e \text{ edge of } \Delta \}.
\end{equation}
Then~$T_v(x) \cong T_{v'}(x)$. In particular, we have~$T_v(x) \cong T_v(y)$ if and only if~$T_{v'}(x) \cong T_{v'}(y)$.
\end{proposition}

This means that in small enough neighbourhoods of vertices, the classification problem of toric fibres does not depend on the choice of vertex. \medskip

\proofof{Proposition~\ref{prop:vertexclass}}
The main idea of the proof is to use the edges of~$\Delta$ to construct symmetric probes exchanging a pair of toric fibres sitting close to the vertices at the endpoints of the given edge. 

Let~$e \subset \Delta$ be an edge of the moment polytope with directional vector~$v_e \in \Lambda^*$. Let~$F,F' \subset \Delta$ be the two facets of~$\Delta$ adjacent to the endpoints of~$e$ but not containing~$e$. We note that the Delzant condition at the vertices adjacent to~$e$ implies that~$v_e$ intersects~$F,F'$ integrally transversely. This can be easily seen by using the corresponding normal form mapping~$e$ to the span of~$e_n^*$ and~$F$ (or~$F'$) to the span of~$e_1^*, \ldots, e_{n-1}^*$. In other words, we obtain symmetric probes parallel to~$e$, as long as both endpoints intersect~$F$ and~$F'$. If~$a < \min\{ \intl(e) \, \vert \, e \text{ edge of } \Delta \}$, then every~$T_v(x) \subset B_v(a)$ can be accessed by such a symmetric probe. Any two vertices~$v,v' \in \Delta$ can be linked by a chain of edges and this proves the claim, up to performing coordinate permutations in one of the two ball embeddings (which can always be realized by symmetric probes). 
\proofend 

The technique used in the previous proof can be generalized to any symmetric probe in a face~$\Delta'$ of a Delzant polytope. Recall that any face of a Delzant polytope is itself Delzant.

\begin{proposition}
Let~$\Delta' \subset \Delta$ be a face and~$\sigma' \subset \Delta'$ a symmetric probe therein. Then there is a neighbourhood~$U$ of~$\sigma'$ such that any parallel translate of~$\sigma'$ in~$U \cap \Int \Delta$ is a symmetric probe. 
\end{proposition}

\proof
Let~$\sigma' \subset \Delta'$ be a symmetric probe with endpoints on the facets~$f,f' \subset \Delta'$. We can write~$f = \Delta' \cap F$ and similarly~$f' = \Delta \cap F'$ for~$F,F'$ facets of~$\Delta$. Let us now show that any parallel translate of~$\sigma'$ with endpoints on~$F$ and~$F'$ is an admissible symmetric probe in~$\Delta$. At the face~$f = \Delta' \cap F$, we can choose a normal form such that~$\Delta'$ spans the coordinate subspace spanned by~$e_1^*,\ldots,e_k^*$ and~$F$ the one spanned by~$e_2^*, \ldots, e_n^*$. Intergral transversality of the intersection of~$\sigma'$ and~$f$ implies that~$v_{\sigma'},e_2^*,\ldots,e_k^*$ is a lattice basis for the sublattice spanned by~$e_1^*,\ldots,e_k^*$. Here, we have denoted the directional vector of the symmetric probe~$\sigma'$ by~$v_{\sigma'}$. This implies that~$v_{\sigma'},e_2^*,\ldots,e_n^*$ is a lattice basis of the full lattice in the ambient space, proving integral transversality.
\proofend

\newpage

\bibliographystyle{abbrv}
\bibliography{mybibfile}

\end{document}